\documentclass{article}

\newtheorem{fed}{\textbf{Definition}}[section]
\newtheorem{thm}[fed]{\textbf{Theorem}}
\newtheorem{lemma}[fed]{\textbf{Lemma}}

\newtheorem{prop}[fed]{\textbf{Proposition}}
\newtheorem{cor}[fed]{\textbf{Corollary}}
\usepackage{amssymb,bbm,graphicx,epsfig,psfrag,epic,eepic,latexsym}
\usepackage{amsmath}
\usepackage{mathrsfs} 
\usepackage[dvips]{color}       

\newcommand{\R}{\mathbb{R}}
\newcommand{\Z}{\mathbb{Z}}
\newcommand{\C}{\mathbb{C}}
\newcommand{\NN}{\mathcal{N}}
\newcommand{\MM}{\mathcal{M}}
\newcommand{\p}{\partial}
\newcommand{\om}{\omega}
\newcommand{\Hess}{{\rm Hess}}
\newcommand{\spec}{{\rm spec}}
\newcommand{\ind}{{\rm ind}}

\begin{document}
\title{A Floer homology for exact contact embeddings}
\author{Kai Cieliebak, Urs Frauenfelder}
\date{4 October 2007}
\maketitle
\tableofcontents
\parindent=0pt
\parskip=4pt

\begin{abstract}
In this paper we construct the Floer homology for an action functional
which was introduced by Rabinowitz and prove a vanishing theorem. As an
application, we show that there 
are no displaceable exact contact embeddings of the unit cotangent
bundle of a sphere of dimension greater than three into a convex exact
symplectic manifold with vanishing first Chern class. This generalizes
Gromov's result that there are 
no exact Lagrangian embeddings of a sphere into $\mathbb{C}^n$. 
\end{abstract}

\tableofcontents

\section[Introduction]{Introduction}

{\bf Exact convex symplectic manifolds and hypersurfaces.}
An {\em exact convex symplectic manifold} $(V,\lambda)$ is a connected
manifold $V$ of dimension $2n$ without boundary with a one-form
$\lambda$ such that the following conditions are satisfied.
\begin{description}
 \item[(i)] The two-form $\omega=d\lambda$ is symplectic.
 \item[(ii)] The symplectic manifold $(V,\omega)$ is convex at
  infinity, i.e.\, there exists an exhaustion $V=\cup_k V_k$ of $V$ by
  compact sets $V_k\subset V_{k+1}$ with smooth boundary such that
  $\lambda|_{\p V_k}$ is a contact form. 
\end{description}
(cf.~\cite{eliashberg-gromov}). Define a vector field $Y_\lambda$ on
$V$ by $i_{Y_\lambda}\om=\lambda$. Then the last condition is
equivalent to saying that $Y_\lambda$ points out of $V_k$ along $\p
V_k$.  

We say that an exact convex symplectic manifold $(V,\lambda)$ is {\em
complete} if the vector field $Y_\lambda$ is complete. We say that
$(V,\lambda)$ has {\em bounded topology} if $Y_\lambda\neq 0$ outside
a compact set. Note that $(V,\lambda)$ is complete and of bounded
topology iff there exists an embedding $\phi:M\times\R_+\to V$ such
that $\phi^*\lambda=e^r\alpha_M$ with contact form
$\alpha_M=\phi^*\lambda|_{M\times\{0\}}$, and such that
$V\setminus\phi(M\times\R_+)$ is compact. (To see this, simply apply
the flow of $Y_\lambda$ to $M:=\p V_k$ for large $k$). 

We say that a subset $A \subset V$ is {\em displaceable} if it can be
  displaced from itself via a Hamiltonian isotopy, i.e.\,there exists
  a smooth family of Hamiltonian functions 
  $H=H(A) \in C^\infty([0,1]\times V)$
  with compact support such that the time one flow $\phi_H$
  of the time dependent Hamiltonian vector field $X_{H_t}$ defined by
  $dH_t=-\iota_{X_{H_t}}\omega$ for $H_t=H(t,\cdot) \in C^\infty(V)$
  and $t \in [0,1]$ satisfies $\phi_H(A) \cap A=\emptyset$. 

The main examples of exact convex symplectic manifolds we have in mind
are Stein manifolds.  We briefly recall its
definition. A {\em Stein manifold} is a triple $(V,J,f)$ where $V$ is
a connected 
manifold, $J$ is an integrable complex structure on $V$ and $f \in
C^\infty(V)$ is an exhausting plurisubharmonic 
function, i.e. $f$ is proper and bounded from below,
and the exact two form $\omega=-d d^c f$ is symplectic. 
Here the
one form $\lambda=-d^c f$ is defined by the condition
$d^c f(\xi)=d f(J \xi)$
for every vector field $\xi$. 
We refer to \cite{cieliebak-eliashberg} for a detailed treatment
of Stein manifolds and Eliashberg's topological characterization of
them. 
It is well known that if the plurisubharmonic function $f$ is
Morse, then all critical points of $f$ have Morse index less than
or equal than half the dimension of $V$, see for example
\cite{cieliebak-eliashberg}. The Stein manifold
$(V,J,f)$ is called \emph{subcritical} if this inequality is strict. 
In a subcritical Stein manifold, every compact subset $A$ is
displaceable~\cite[Lemma 3.2]{biran-cieliebak}. 

{\em Remark. }
Examples of exact convex symplectic manifolds which are not Stein can
be obtained using the following construction. Let $M$ be a
$(2n-1)$-dimensional closed manifold which admits a pair of contact forms
$(\alpha_0,\alpha_1)$ satisfying
$$\alpha_1 \wedge (d\alpha_1)^{n-1}=-\alpha_0 \wedge (d\alpha_0)^{n-1}>0$$
and
$$\alpha_i\wedge(d\alpha_i)^k \wedge(d\alpha_j)^{n-k-1}=0, \quad
0 \leq k \leq n-2$$
where $(i,j)$ is a permutation of $(0,1)$. Then a suitable interpolation
of $\alpha_0$ and $\alpha_1$ endows the manifold $V=M \times [0,1]$
with the structure of an exact convex symplectic manifold, where
the restriction of the one-form to $M \times \{0\}$ is given by
$\alpha_0$ and the restriction to $M \times \{1\}$ is given by
$\alpha_1$. Since $H_{2n-1}(V)=\Z$, the manifold $V$ does not admit a
Stein structure. Moreover, what
makes these examples particularly interesting is the fact that
they have two boundary components, whereas the boundary of a connected
Stein manifold is always connected. The first construction in
dimension four of an exact convex
symplectic manifold of the type above was carried out by D.\,McDuff
in \cite{mcduff}. H.\,Geiges generalized her method in \cite{geiges},
where he also obtained higher dimensional examples. 

If $(V,\lambda)$ is an exact convex symplectic manifold then so is
its {\em stabilization} $(V \times \mathbb{C}, \lambda \oplus
\lambda_\C)$ for the one form $\lambda_\C=\frac{1}{2}(x\,dy-y\,dx)$ on
$\mathbb{C}$. Moreover, in $(V \times \mathbb{C}, \lambda \oplus
\lambda_\C)$ every compact subset $A$ is displaceable. 
It is shown in \cite{cieliebak} that each subcritical Stein manifold
is Stein deformation equivalent to a split Stein manifold, i.e.\,a
Stein manifold of the form $(V \times \mathbb{C},J \times i,
f+|z|^2)$ for a Stein manifold $(V,J,f)$.

{\em Remark. }If $(V,\lambda)$ is an exact convex symplectic manifold,
then so is $(V,\lambda+dh)$ for any smooth function $h:V\to\R$ with compact
support. We call the 1-forms $\lambda$ and $\lambda+dh$ {\em
  equivalent}. For all our considerations only the equivalence class of
$\lambda$ will be relevant. 
\medskip

An {\em exact convex hypersurface} in an exact convex symplectic
manifold $(V,\lambda)$ is a compact hypersurface (without boundary)
$\Sigma\subset V$ such that
\begin{description}
\item [(i)]There exists a contact 1-form $\alpha$ on $\Sigma$ such that
$\alpha-\lambda|_\Sigma$ is exact.
\item[(ii)] $\Sigma$ is {\em bounding},  
i.e. $V\setminus\Sigma$ consists of two connected components, one
compact and one noncompact. 
\end{description}

{\em Remarks. }
(1) It follows that the volume form $\alpha\wedge(d\alpha)^{n-1}$
defines the orientation of $\Sigma$ as boundary of the bounded
component of $V\setminus\Sigma$. 

(2) If $\Sigma$ is an exact convex hypersurface in $(V,\lambda)$ with
contact form $\alpha$, then there exists an equivalent 1-form
$\mu=\lambda+dh$ on $V$ such that $\alpha=\mu|_\Sigma$. To see this,
extend $\alpha$ to a 1-form $\beta$ on $V$. As
$(\beta-\lambda)|_\Sigma$ is exact, there exists a function $h$ on a
neighbourhood $U$ of $\Sigma$ such that $\beta-\lambda=dh$ on $U$. Now
simply extend $h$ to a function with compact support on $V$ and set
$\mu:=\lambda+dh$. 

(3) If $H^1(\Sigma;\R)=0$ condition (i) is equivalent to
$d\alpha=\om|_\Sigma$. 

(4) Condition (ii) is automatically satisfied if $H_{2n-1}(V;\Z)=0$,
e.g.~if $V$ is a stabilization or a Stein manifold of dimension $>2$.  
\medskip

{\bf Floer homology. }
In the following we assume that $(V,\lambda)$ is a complete exact
convex symplectic manifold of bounded topology, and $\Sigma\subset V$
is an exact convex hypersurface with contact form $\alpha$. 
We will define an invariant $HF(\Sigma,V)$ as the Floer homology of an
action functional which was studied previously by Rabinowitz
\cite{rabinowitz}.

A {\em defining Hamiltonian for $\Sigma$} is a function $H \in C^\infty(V)$  
which is  constant outside of a compact set of $V$,
whose zero level set $H^{-1}(0)$ equals $\Sigma$,
and whose Hamiltonian vector field
$X_H$ defined by $dH=-\iota_{X_H}\omega$ agrees with the Reeb vector field
$R$ of $\alpha$ on $\Sigma$. Defining Hamiltonians exist since $\Sigma$ is
bounding, and they form a convex space.  

Fix a defining Hamiltonian $H$ and 
denote by $\mathscr{L}=C^\infty(\mathbb{R}/\mathbb{Z},V)$ 
the free loop space of $V$. 
Rabinowitz' action functional 
$$\mathcal{A}^H \colon \mathscr{L}\times \mathbb{R} \to \mathbb{R}$$
is defined as
$$\mathcal{A}^H(v,\eta) := \int_0^1v^*\lambda-\eta
\int_0^1 H(v(t))dt,\quad (v,\eta) \in \mathscr{L}\times \mathbb{R}.$$
One may think of $\mathcal{A}^H$ as the Lagrange multiplier functional
of the unperturbed action functional of classical mechanics also studied
in Floer theory to a mean value constraint of the loop. The critical
points of $\mathcal{A}^H$ satisfy 
\begin{equation}\label{crit1}
\left. \begin{array}{cc}
\partial_t v(t)= \eta X_H(v(t)),
& t \in \mathbb{R}/\mathbb{Z}, \\
H(v(t))=0. & \\
\end{array}
\right\}
\end{equation}
Here we used the fact that $H$ is invariant under its Hamiltonian
flow. Since the restriction of the Hamiltonian vector
field $X_H$ to $\Sigma$ is the Reeb vector field, 
the equations (\ref{crit1}) are equivalent to
\begin{equation}\label{crit}
\left. \begin{array}{cc}
\partial_t v(t)= \eta R(v(t)),
& t \in \mathbb{R}/\mathbb{Z}, \\
v(t)\in \Sigma, & t \in \mathbb{R}/\mathbb{Z}, \\
\end{array}
\right\}
\end{equation}
i.e. $v$ is a periodic orbit of the Reeb vector field on
$\Sigma$ with period $\eta$.
\footnote{The period $\eta$ may be negative or zero. 
We refer in this paper
to Reeb orbits moved backwards as Reeb orbits with negative period
and to constant orbits as Reeb orbits of period zero.}


\begin{thm}\label{well}
Under the above hypotheses, the Floer homology
$HF(\mathcal{A}^H)$ is well-defined. Moreover, if 
$H_s$ for $0 \leq s \leq 1$ is a smooth family of defining functions
for exact convex hypersurfaces $\Sigma_s$, then
$HF(\mathcal{A}^{H_0})$ and $HF(\mathcal{A}^{H_1})$ are
canonically isomorphic.
\end{thm}

Hence the Floer homology $HF(\mathcal{A}^H)$ is independent of
the choice of the defining function $H$ for an exact convex
hypersurface $\Sigma$, and the resulting invariant
$$
   HF(\Sigma,V) := HF(\mathcal{A}^H)
$$
does not change under homotopies of exact convex hypersurfaces. 

The next result is a vanishing theorem for the Floer homology
$HF(\Sigma,V)$.

\begin{thm}\label{zero}
If $\Sigma$ is displaceable, then $HF(\Sigma,V)=0$.
\end{thm}

{\em Remark. }The action functional $\mathcal{A}^H$ is also defined if
$H^{-1}(0)$ is not exact convex. However, in this case the Floer homology
$HF(\mathcal{A}^H)$ cannot in general be defined because the
moduli spaces of flow lines will in general not be compact
up to breaking anymore. The problem is that the Lagrange multiplier $\eta$
may go to infinity. This phenomenon actually does happen as the
counterexamples to the Hamiltonian Seifert conjecture show, see
\cite{ginzburg-gurel} and the literature cited therein.
\medskip

Denote by $c_1$ the first Chern class of the tangent bundle of $V$
(with respect to an $\om$-compatible almost complex structure and
independent of this choice, see~\cite{mcduff-salamon}). Evaluation of
$c_1$ on spheres gives rise to a 
homomorphism $I_{c_1} \colon \pi_2(V) \to \mathbb{Z}$. 
If $I_{c_1}$ vanishes then the Floer homology $HF_*(\Sigma,V)$ can be
$\mathbb{Z}$-graded with half integer degrees, i.e. $* \in 1/2+\mathbb{Z}$.

The third result is a computation of the Floer homology for the unit
cotangent bundle of a sphere. 

\begin{thm}\label{compute}
Let $(V,\lambda)$ be a complete exact convex symplectic
manifold of bounded topology satisfying $I_{c_1}=0$. Suppose that
$\Sigma\subset V$ is an exact convex hypersurface with contact form
$\alpha$ such that $(\Sigma,\ker\alpha)$ is contactomorphic to the
unit cotangent bundle $S^*S^n$ of the sphere of dimension $n\geq 4$
with its standard contact structure. Then  
$$
HF_k(\Sigma,V) = \left\{ \begin{array}{cc}
\mathbb{Z}_2,
& k \in \{-n+\frac{1}{2},-\frac{1}{2},\frac{1}{2},n-\frac{1}{2}\}+
\mathbb{Z}\cdot(2n-2), \\
0, & \mathrm{else.} \\
\end{array}
\right.
$$
\end{thm}
\medskip

{\bf Applications and discussion. }
The following well-known technical lemma will allow us to remove
completeness and bounded topology from the hypotheses of our
corollaries.

\begin{lemma}\label{bt}
Assume that $\Sigma$ is an exact convex hypersurface in the exact convex
symplectic manifold $(V,\lambda)$. Then $V$ can be modified outside of
$\Sigma$ to an exact convex symplectic manifold
$(\hat{V},\hat\lambda)$ which is complete and of 
bounded topology. If $I_{c_1}=0$ for $V$ the same holds for
$\hat V$. If $\Sigma$ is displaceable in $V$, then we can
arrange that it is displaceable in $\hat V$ as well. 
\end{lemma}

\textbf{Proof: }
Let $V_1\subset V_2\dots$ be the compact exhaustion 
in the definition of an exact convex symplectic manifold. 
Since $\Sigma$ is compact, it is contained in $V_k$ for some $k$. The
flow of $Y_\lambda$ for times $r\in(-1,0]$ defines an embedding $\phi:\p
V_k\times(-1,0]\to V_k$ such that $\phi^*\lambda=e^r\lambda_0$,
where $\lambda_0=\lambda|_{\p V_k}$. Now define
$$
   (\hat V,\hat\lambda) := (V_k,\lambda)\cup_\phi\bigl(\p
   V_k\times(-1,\infty),e^r\lambda_0\bigr). 
$$
This is clearly complete and of bounded topology. The statement about
$I_{c_1}$ is obvious. If $\Sigma$
is displaceable by a Hamiltonian isotopy generated by a compactly
supported Hamiltonian $H:[0,1]\times V\to\R$, we choose $k$ so large
that ${\rm supp}H\subset[0,1]\times V_k$ and apply the same
construction.  
\hfill $\square$

As a first consequence of Theorem~\ref{zero}, we recover some known 
cases of the Weinstein conjecture, see~\cite{viterbo},
\cite{frauenfelder-ginzburg-schlenk}. 

\begin{cor}
Every displaceable exact convex hypersurface $\Sigma$ in an exact
convex symplectic manifold $(V,\lambda)$ carries a closed
characteristic. In particular, 
this applies to all exact convex hypersurfaces in a subcritical Stein
manifold, or more generally in a stabilization $V\times\C$. 
\end{cor}

\textbf{Proof: }
In view of Lemma~\ref{bt}, we may assume without loss of generality
that $(V,\lambda)$ is complete and of bounded topology. Then by
Theorem~\ref{zero} the Floer homology $HF(\mathcal{A}^H)$ vanishes,
where $H$ is a defining function for $\Sigma$. On the other hand, 
the action functional $\mathcal{A}^H$ always has critical points
corresponding to the constant loops in $\Sigma$. So the vanishing of
the Floer homology implies that there must also exist nontrivial
solutions of (\ref{crit}), which are just closed characteristics, 
connected to constant loops by gradient flow lines of $\mathcal{A}^H$.  
\hfill $\square$

For further applications, the following notation will be convenient. 
An {\em exact contact embedding} of a closed contact manifold
$(\Sigma,\xi)$ into an exact convex symplectic manifold $(V,\lambda)$
is an embedding $\iota \colon \Sigma \to V$ such that
\begin{description}
\item [(i)]There exists a 1-form $\alpha$ on $\Sigma$ such that
  $\ker\alpha=\xi$ and $\alpha-\iota^*\lambda$ is exact.
\item[(ii)] The image $\iota(\Sigma)\subset V$ is bounding. 
\end{description}
In other words, $\iota(\Sigma)\subset V$ is an exact convex
hypersurface with contact form $\iota_*\alpha$ which is
contactomorphic (via $\iota^{-1}$) to $(\Sigma,\xi)$.

Now Theorems~\ref{zero} and~\ref{compute} together with Lemma~\ref{bt}
immediately imply

\begin{cor}\label{cor:non-displaceable}
Assume that $n\geq 4$ and
there exists an exact contact embedding $\iota$ of
$S^*S^n$ into an exact convex symplectic manifold satisfying
$I_{c_1}=0$. Then $\iota(S^*S^n)$ is not displaceable.
\end{cor}

Since in a stabilization $V\times\C$ all compact subsets are
displaceable, we obtain in particular

\begin{cor}\label{cor:subcrit}
For $n\geq 4$ there does not exist an exact contact embedding of $S^*S^n$
into a subcritical Stein manifold, or more generally, into the
stabilization $(V\times\C,\lambda\oplus\lambda_\C)$ of an exact convex
symplectic manifold $(V,\lambda)$ satisfying $I_{c_1}=0$.  
\end{cor}

{\em Remark. }
If $n$ is even then there are no smooth embeddings of
$S^*S^n$ into a subcritical Stein manifold by topological reasons, see
Appendix~\ref{app:top}. However, at least for $n=3$ and $n=7$ there
are no topological obstructions, see the discussion below. 
\medskip

If $(V,J,f)$ is a Stein manifold with $f$ a Morse function,
P.~Biran~\cite{biran} defines the {\em critical coskeleton} as the union
of the unstable manifolds (w.r.~to $\nabla f$) of the critical points
of index $\dim V/2$. It is proved in~\cite{biran} that every compact
subset $A\subset V$ which does not intersect the critical coskeleton
is displaceable. For example, in a cotangent bundle the critical
coskeleton (after a small perturbation) is one given fibre. Thus 
Corollary~\ref{cor:non-displaceable} implies

\begin{cor}\label{cor:crit}
Assume that $n\geq 4$ and there exists an exact contact
embedding $\iota$ of $S^*S^n$ into a Stein manifold $(V,J,f)$
satisfying $I_{c_1}=0$. Then
$\iota(\Sigma)$ must intersect the critical coskeleton. In particular,
the image of an exact contact embedding of $S^*S^n$ into a cotangent
bundle $T^*Q$ must intersect every fibre. 
\end{cor}

{\em Remark. }
Let $\iota \colon L \to V$ be an {\em exact Lagrangian embedding} of
$L$ into $V$, i.e. such that 
$\iota^*\lambda$ is exact. Since by
Weinsteins Lagrangian neighbourhood theorem 
\cite[Theorem 3.33]{mcduff-salamon1} a tubular neighbourhood of
$\iota(L)$ can be symplectically identified with a tubular neighbourhood
of the zero section of the cotangent bundle of $L$, we obtain
an exact contact embedding of $S^*L$ into $V$. Thus the last 3
corollaries generalize
corresponding results about exact Lagrangian embeddings. For example,
Corollary~\ref{cor:subcrit} generalizes (for spheres) the well-known
result~\cite{gromov,audin-lalonde-polterovich,biran-cieliebak} that
there exist no exact Lagrangian embeddings into subcritical Stein
manifolds. Corollary~\ref{cor:crit} implies (cf.~\cite{biran}) that
an embedded Lagrangian sphere of dimension $\geq 4$ in a cotangent
bundle $T^*Q$ must intersect every fibre.
\medskip

{\em Remark. }
Let us discuss Corollary~\ref{cor:subcrit} in the cases $n\leq 3$ that
are not accessible by our method of proof. We always equip $\C^n$ with
the canonical 1-form
$\lambda=\frac{1}{2}\sum_{i=1}^n(x_idy_i-y_idx_i)$. 

$n=1$: Any embedding of two disjoint circles into $\C$ is an exact
contact embedding of $S^*S^1$, so Corollary~\ref{cor:subcrit} fails in
this case. 

$n=2$: In this case Corollary~\ref{cor:subcrit} is true for purely
topological reasons; we present various proofs in
Appendix~\ref{app:top}. 

$n=3$: In this case Corollary~\ref{cor:subcrit} is true for
subcritical Stein manifolds and can be
proved using symplectic homology, see the last remark in this section.  
\medskip

{\em Example. }
In this example we illustrate that the preceding results about exact
contact embeddings are sensitive to the contact structure. Let $n=3$
or $n=7$. Then $S^*S^n\cong S^n\times
S^{n-1}$ embeds into $\mathbb{R}^{n+1} \times S^{n-1}$. 
On the other hand $\mathbb{R}^{n+1} \times S^{n-1}$ is 
diffeomorphic to
the subcritical Stein manifold $T^*S^{n-1}\times\C$,
and identifying $S^*S^n$ with a level set in
$T^*S^{n-1}\times \C$ defines a contact structure $\xi$ on $S^*S^n$ . Thus
$(S^*S^n,\xi)$ has an exact contact embedding into a subcritical Stein
manifold (in fact into $\C^n$) for $n=3,7$, whereas $(S^*S^7,\xi_{\rm
  st})$ admits no such embedding by Corollary~\ref{cor:subcrit}. In
particular, we conclude 

\begin{cor}\label{cor:non-diffeo}
The two contact structures $\xi$ and $\xi_{\rm st}$ on $S^*S^7\cong
S^7\times S^6$ described above are not diffeomorphic.
\end{cor}

{\em Remarks. }
(1) Corollary~\ref{cor:non-diffeo} also holds in the case $n=3$, although
our method does not apply there. Indeed, the contact structures $\xi$
and $\xi_{\rm st}$ on $S^*S^n$ for $n=3,7$ are distinguished by their
cylindrical contact homology (see~\cite{ustilovsky}, \cite{yau}). 

(2) The contact structures $\xi$ and $\xi_{\rm st}$ on $S^3\times S^2$
are homotopic as almost contact structures, i.e.~as symplectic
hyperplane distributions. This follows simply from the fact
(see e.g.~\cite{geiges-book}) that on 5-manifolds almost contact
structures are classified up to homotopy by their first Chern classes
and $c_1(\xi)=c_1(\xi_{\rm st})=0$. It would be interesting to know
whether $\xi$ and $\xi_{\rm st}$ on $S^7\times S^6$ are also homotopic
as almost contact structures. Here the first obstruction to such a
homotopy vanishes because $c_3(\xi)=c_3(\xi_{\rm st})=0$, but there
are further obstructions in dimensions $7$ and $13$ which remain to be
analysed along the lines of~\cite{morita}. 


{\em Remark }(obstructions from symplectic field theory). 
Symplectic field theory~\cite{eliashberg-givental-hofer} also yields
obstructions to exact 
contact embeddings. For example, by neck stretching along the image of
an exact contact embedding, the following result is proved
in~\cite{cieliebak-mohnke}: {\em Let $(\Sigma^{2n-1},\xi)$ be a closed
contact manifold with $H_1(\Sigma;\Z)=0$ which admits
an exact contact embedding into $\C^n$. Then for every nondegenerate
contact form defining $\xi$ there exist closed Reeb orbits of
Conley-Zehnder indices $n+1+2k$ for all integers $k\geq 0$.}\\
Here Conley-Zehnder indices are defined with respect to
trivializations extending over spanning surfaces. This result applies
in particular to the unit cotangent bundle $\Sigma=S^*Q$ of a closed
Riemannian manifold $Q$ with $H_1(Q;\Z)=0$. For example, if $Q$
carries a metric of nonpositive curvature then all indices are $\leq
n-1$ and hence $S^*Q$ admits no exact contact embedding into
$\C^n$. On the other hand, any nondegenerate metric on the sphere
$S^n$ has closed geodesics of all indices $n+1+2k$, $k\geq 0$, so this
result does {\em not} exclude exact contact embeddings
$S^*S^n\hookrightarrow\C^n$.  

{\em Remark }(obstructions from symplectic homology). 
Corollary~\ref{cor:subcrit} for subcritical Stein manifolds can be
proved for all $n\geq 3$ by combining the following five results on
symplectic homology. See~\cite{cieliebak-oancea} for details.

(1) The symplectic homology $SH(V)$ of a subcritical Stein manifold
$V$ vanishes~\cite{cieliebak-chord}. 

(2) If $\Sigma\subset V$ is an exact convex hypersurface in an
exact convex symplectic manifold bounding the compact domain $W\subset
V$, then $SH(V)=0$ implies $SH(W)=0$. This follows from an argument by
M.~McLean~\cite{mclean}, based on Viterbo's transfer map~\cite{viterbo}
and the ring structure on symplectic homology. 

(3) If $SH(W)=0$, then the positive action part $SH^+(W)$ of
symplectic homology is only nonzero in finitely many degres. This
follows from the long exact sequence induced by the action
filtration. 

(4) $SH^+(W)$ equals the non-equivariant linearized contact homology
$NCH(W)$. This is implicit in~\cite{bourgeois-oancea}, see
also~\cite{cieliebak-oancea}. 

(5) If $\p W=S^*S^n$ and $n\geq 3$, then $NCH(W)$ is independent of
the exact filling $W$ and equals the homology of the free loop space
of $S^n$ (modulo the constant loops), which is nonzero in infinitely
many degrees.

\section{Exact contact embeddings}\label{sec:exact}

Let $\Sigma$ be a connected 
closed $2n-1$ dimensional manifold. A {\em contact structure}
$\xi$ is a field of hyperplanes $\xi \subset T \Sigma$ such that
there exists a one-form $\alpha$ satisfying
$$\xi = \ker \alpha, \quad \alpha \wedge d\alpha^{n-1} >0.$$ 
The one form $\alpha$ is called a contact form. It is determined by 
$\xi$ up to multiplication with a function $f>0$.  Given a contact
form $\alpha$ the {\em Reeb vector field} $R$ on $\Sigma$ is defined by
the conditions
$$\iota_R d\alpha=0, \quad \alpha(R)=1.$$
Unit cotangent bundles have a natural contact structure as the following
example shows.

{\em Example. }
For a manifold $N$ we denote by $S^*N$ the oriented projectivization of
its cotangent bundle $T^*N$, i.e. elements of $S^*N$ are 
equivalence classes $[v^*]$ of contangent vectors $v^* \in T^* N$ under
the equivalence relation $v^* \cong w^*$ if there
exists $r>0$ such that $v^*=r w^*$. Denote by $\pi \colon S^* N \to N$
the canonical projection. A contact structure $\xi$ on $S^*N$ is
given by
$$\xi_{[v^*]}=\ker v^* \circ d\pi([v^*]).$$
If $g$ is a Riemannian metric on $N$ then
$S^*N$ can be identified with the space of tangent vectors
of $N$ of length one and the restriction of the Liouville one form
defines a contact form. Observe that the Reeb vector field generates 
the geodesic flow on $N$.

If $\iota \colon \Sigma \to V$ is a exact contact embedding, then
$\alpha=\iota^* \lambda$ defines a contact form for the
contact structure $\xi$. One might ask which contact forms $\alpha$
can arise in this way. The following proposition shows that if one
contact form defining the contact structure $\xi$ arises from an exact
contact embedding, then every other contact form defining $\xi$ arises
as well. 

\begin{prop} \label{embe}
Assume that $\iota:(\Sigma, \xi)\to (V,\lambda)$ is an exact contact
embedding with $\xi=\ker\iota^*\lambda$. Then for every contact form
$\alpha$ defining the contact structure $\xi$ on $\Sigma$ there exists
a constant $c>0$ and a bounding  embedding $\iota_\alpha \colon \Sigma
\to V$ such that $\iota_\alpha^* \lambda=c\alpha$. 
\end{prop} 
\textbf{Proof of Proposition~\ref{embe}: } 
The proof uses the fact that if there exists
an exact contact embedding of a contact manifold 
into an exact convex symplectic manifold $(V,\lambda)$ then the
negative symplectization can be embedded. To see that we need two
facts. Recall that the vector field $Y_\lambda$ is defined by the
condition $\lambda=\iota_{Y_\lambda}d \lambda$. 

{\em Fact 1: The flow $\phi_\lambda^t$ of $Y_\lambda$
exists for all negative times $t$.}
\\
Indeed, let $x\in V$. Then $x\in
V_k$ for some $k$. As $Y_\lambda$ points out of $V_k$ along $\p V_k$,
$\phi_\lambda^t(x)\in V_k$ for all $t\leq 0$ and compactness of $V_k$
yields completeness for $t\leq 0$.
\\ \\
{\em Fact 2: The vector field $Y_\lambda$ satisfies 
\begin{equation}\label{liouville}
\iota_{Y_\lambda}\lambda=0, \quad L_{Y_\lambda}\lambda=\lambda,
\end{equation}
where $L_{Y_\lambda}$ is the Lie derivative along the vector field
$Y_\lambda$. In particular, the flow $\phi^r_\lambda$ of $Y_\lambda$ 
satisfies $(\phi^r_\lambda)^*\lambda=e^r \lambda$.}
\\
Indeed, the first equation in (\ref{liouville}) follows directly
from the definition of $Y_\lambda$. To prove the second one
we compute using Cartan's formula
$$L_{Y_\lambda}\lambda=d\iota_{Y_\lambda}\lambda+\iota_{Y_\lambda}d\lambda
=\lambda.$$
\\
Now set $\alpha_0=\iota^* \lambda$ and 
consider the symplectic manifold 
$\big(\Sigma \times \mathbb{R}_-, d (e^r \alpha_0)\big)$ where $r$
denotes the coordinate on the $\mathbb{R}$-factor. By Fact 1, the flow
$\phi_\lambda^r$ exists for all $r\leq 0$. By Fact 2, the embedding 
$$\hat{\iota} \colon \Sigma \times \mathbb{R}_- \to V, \quad
(x,r) \mapsto \phi_\lambda^r(\iota(x))$$
satisfies
$$(\hat{\iota})^* \lambda=e^r \alpha_0.$$
If $\alpha$ is another
contact form on $\Sigma$ which defines the contact structure
$\xi$ then there exists a smooth function 
$\rho_\alpha \in C^\infty(\Sigma)$ such that
$$\alpha=e^{\rho_\alpha} \alpha_0.$$
Set $m:=\max_\Sigma\rho_\alpha$ and $c:=e^{-m}$. Then 
$$\iota_\alpha \colon \Sigma \to V, \quad
x \mapsto \hat{\iota}(x,\rho_\alpha(x)-m)$$
gives the required contact embedding for $\alpha$.
This proves the proposition. \hfill $\square$

{\em Remark. }
If the vector field $Y_\lambda$ is complete, then the preceding proof
yields a symplectic embedding of the whole symplectization
$\bigl(\Sigma\times\R,d(e^r\alpha_0)\bigr)$ into $(V,\om)$.

\section{Floer homology for Rabinowitz's action
functional}\label{sec:floer} 

In this section we construct the Floer homology for Rabinowitz's action
functional defined in the introduction and prove Theorem~\ref{well}
and Theorem~\ref{zero}. We assume that the reader is familiar
with the constructions in Floer theory which can be found in
Floer's original papers \cite{floer1,floer2, floer3, floer4, floer5} or
in Salamon's lectures \cite{salamon}. The finite dimensional case
of Morse theory is treated in the book of Schwarz \cite{schwarz}. 

Throughout this section we maintain the following setup:
\begin{itemize}
\item $(V,\lambda)$ is a complete exact convex symplectic manifold of
  bounded topology.
\item $\Sigma\subset V$ is an exact convex hypersurface with contact
  form $\alpha$ and defining Hamiltonian $H$.
\end{itemize}
Our sign conventions for Floer homology are as follows:
\begin{itemize}
\item The {\em Hamiltonian vector field} $X_H$ is defined by
$dH=-i_{X_H}\om$, where $\om=d\lambda$. 
\item An almost complex structure $J$ on $V$ is {$\om$-compatible} if
  $\om(\cdot,J\cdot)$ defines a Riemannian metric. Thus
  the gradient with respect to this metric is related to the
  symplectic vector field by $X_H=J\nabla H$. 
\item Floer homology is defined using the {\em positive} gradient flow
  of the action functional $\mathcal{A}^H$.  
\end{itemize}

The action functional $\mathcal{A}^H$ is invariant under
the $S^1$-action on $\mathscr{L}\times \mathbb{R}$ given
by $t_*(v(\cdot),\eta) \mapsto (v(t+\cdot),\eta)$. In particular,
the action functional $\mathcal{A}^H$ will never be Morse.
However, generically it is {\em Morse-Bott}, i.e.
its critical set is a manifold whose tangent space is the kernel of
the Hessian of the action functional. 
We make the following nondegeneracy assumption on the Reeb flow
$\phi_t$ of the contact form $\alpha$ on $\Sigma$.
\begin{itemize}
 \item[(A)] The closed Reeb orbits of $(\Sigma,\alpha)$ are of {\em
  Morse-Bott type}, i.e.~for each $T\in\R$ the set
  $\NN_T\subset\Sigma$ formed by the $T$-periodic Reeb orbits is a
  closed submanifold, the rank of $d\alpha|_{\NN_T}$ is locally
  constant, and $T_p\NN_T=\ker(T_p\phi_T-{\rm id})$ for all $p\in\NN_T$.     
\end{itemize}
If the assumption (A) does not hold we consider a hypersurface
close by. Note that the contact condition is an open condition and
the assumption (A) is generically satisfied. Since we prove that
our homology is invariant under homotopies we can assume without
loss of generality that (A) holds.  
If (A) is satisfied, then the action functional $\mathcal{A}^H$
is Morse-Bott. 

{\em Remark. }Generically, we can even achieve that all $T$-periodic
Reeb orbits $\gamma$ with $T\neq 0$ are {\em nondegenerate},
i.e.~the linearization $T_p\phi_T:\xi_p\to\xi_p$ at $p=\gamma(0)$ does
not have 1 in its spectrum. In this case the critical manifold of
$\mathcal{A}^H$ consists of a union of
circles for each nonconstant Reeb orbit and a copy of 
the hypersurface $\Sigma$ for the constant solutions, i.e.\,critical 
points with $\eta=0$. Moreover, observe that a nonconstant
Reeb orbit gives rise to infinitely many of them because
it can be repeatedly passed and also be passed in the backward direction.

There are several ways to deal with Morse-Bott situations in Floer
homology. One possibility is to choose an additional small perturbation
to get a Morse situation. This was carried out by Pozniak \cite{pozniak},
where it was also shown that the local Floer homology near each
critical manifold coincides with the Morse homology of the critical 
manifold. Another possibility is to choose an additional Morse function
on the critical manifold. The chain complex is then generated by
the critical points of this Morse function while the boundary operator
is defined by counting flow lines with cascades. This approach
was carried out by the second named author in \cite{frauenfelder}.

Cascades are finite energy gradient flow lines of the action functional
$\mathcal{A}^H$. In the Morse-Bott case the finite energy
assumption is equivalent to assume that the gradient flow line
converges at both ends exponentially to a point on the critical manifold. 
In order to prove that the Floer homology is well defined one
has to show that the moduli spaces of cascades are compact modulo breaking.
There are three difficulties one has to solve.
\begin{itemize}
 \item An $L^\infty$-bound on the loop $v \in \mathscr{L}$.
 \item An $L^\infty$-bound on the Lagrange multiplier $\eta \in \mathbb{R}$.
 \item An $L^\infty$-bound on the derivatives of the loop $v$.
\end{itemize}
Although the first and the third point are nontrivial 
they are standard problems in Floer theory one knows how to deal with. The 
$L^\infty$-bound for the loop follows from the convexity assumption on
$V$ and the derivatives can be controlled since 
our symplectic manifold  is exact and hence
there is no bubbling of pseudo-holomorphic spheres. 
The new feature is the bound on
the Lagrange multiplier $\eta$. We will explain in detail how this can be
achieved. It will be essential that our hypersurface is convex. 
\\ \\
We first explain the metric on the space $\mathscr{L}\times \mathbb{R}$
and deduce from that the equation for the cascades. 
The metric on $\mathscr{L}\times \mathbb{R}$ is the product metric
of the standard metric on $\mathbb{R}$ and a metric on
$\mathscr{L}$ coming from a family of $\omega$-compatible almost
complex structures $J_t$ on $V$. 
For such a family of $\omega$-compatible almost complex structures
$J_t$ we define the metric $g_J$ on $\mathscr{L}\times \mathbb{R}$
for $(v,\eta) \in \mathscr{L}\times \mathbb{R}$ and
$(\hat{v}_1,\hat{\eta}_1), (\hat{v}_2,\hat{\eta}_2) \in
T_{(v,\eta)}(\mathscr{L}\times \mathbb{R})$ by 
$$g_J\big((\hat{v}_1,\hat{\eta}_1),(\hat{v}_2,\hat{\eta}_2)\big)
=\int_0^1 \omega(\hat{v}_1,J_t(v)\hat{v}_2) dt +
\hat{\eta}_1 \cdot \hat{\eta}_2.$$
The gradient of $\mathcal{A}^H$ with respect to this metric is given by
$$
   \nabla \mathcal{A}^H=\nabla_{g_J} \mathcal{A}^H =
   \left(\begin{array}{cc}
   -J_t(v)\bigl(\partial_t v-\eta X_H(v)\bigr)\\
   -\int_0^1 H(v(\cdot,t) dt
   \end{array}\right)\;.
$$
Thus gradient flow lines of $\nabla \mathcal{A}^H$ are solutions 
$(v,\eta) \in C^\infty(\mathbb{R}\times S^1,V \times \mathbb{R})$
of the following problem
\begin{equation}\label{flowline}
\left. \begin{array}{cc}
\partial_s v+J_t(v)(\partial_t v-\eta X_H(v))=0\\
\partial_s \eta+\int_0^1 H(v(\cdot,t)dt=0.
\end{array}
\right\}
\end{equation}
The following proposition is our main tool to bound the Lagrange
multiplier $\eta$.

\begin{prop}\label{mainest}
There exists $\epsilon>0$ such that for every $M > 0$
there exists a constant $c_M<\infty$ such that
$$
\left\{ \begin{array}{c}
||\nabla\mathcal{A}^H(v,\eta)||\leq\epsilon \\
|\mathcal{A}^H(v,\eta)|\leq M \\
\end{array}
\right\}
\quad \Longrightarrow \quad |\eta| \leq c_M.
$$
\end{prop}

We first prove a lemma which says that the action value of a critical
point of $\mathcal{A}^H$, i.e.\,a Reeb orbit,
is given by the period.

\begin{lemma}\label{period}
Let $(v,\eta) \in \mathrm{crit}(\mathcal{A}^H)$, then
$\mathcal{A}^H(v,\eta)=\eta$.
\end{lemma} 

\textbf{Proof: } Inserting (\ref{crit}) into $\mathcal{A}^H$ we
compute
$$\mathcal{A}^H(v,\eta)=\eta\int_0^1\lambda(v)R(v)=\eta
\int_0^1\alpha(v)R(v)=\eta.$$
This proves the lemma. \hfill $\square$
\\

\textbf{Proof of Proposition~\ref{mainest}: }
We prove the proposition in three steps. The first step is an
elaboration of the observation in Lemma~\ref{period}.

\textbf{Step 1: }\emph{There exists $\delta>0$ and a constant
$c_\delta<\infty$ with the following property. For every
$(v,\eta) \in \mathscr{L}\times \mathbb{R}$ such that
$v(t) \in U_\delta=H^{-1}\big((-\delta,\delta)\big)$ for every 
$t \in \mathbb{R}/\mathbb{Z}$, the following estimate holds:
$$|\eta| \leq 2|\mathcal{A}^H(v,\eta)|+c_\delta||
\nabla\mathcal{A}^H(v,\eta)||.$$
}

Choose $\delta>0$ so small such that
$$\lambda(x)X_H(x) \geq \frac{1}{2}+\delta, 
\quad x \in U_\delta.$$
Set
$$c_\delta=2||\lambda|_{U_\delta}||_\infty.$$
We estimate
\begin{eqnarray*}
|\mathcal{A}^H(v,\eta)|&=&
\Bigg|\int_0^1\lambda(v)\partial_t v-\eta\int_0^1 H(v(t))dt\Bigg|\\
&=&\Bigg| \eta \int_0^1\lambda(v)X_H(v)
+\int_0^1\lambda(v)\big(\partial_t v-\eta X_H(v)\big)
-\eta\int_0^1 H(v(t))dt\Bigg| \\
&\geq&
\Bigg| \eta \int_0^1\lambda(v)X_H(v)\Bigg|
-\Bigg|\int_0^1\lambda(v)\big(\partial_t v-\eta X_H(v)\big)\Bigg|\\
& &-\Bigg|\eta\int_0^1 H(v(t))dt\Bigg| \\
&\geq& |\eta|\bigg(\frac{1}{2}+\delta\bigg)-
\frac{c_\delta}{2}||\partial_t v-\eta X_H(v)||_1
-|\eta|\delta\\
&\geq& \frac{|\eta|}{2}-\frac{c_\delta}{2}||\partial_t v-\eta X_H(v)||_2\\
&\geq& \frac{|\eta|}{2}-\frac{c_\delta}{2}||\nabla \mathcal{A}^H(v,\eta)||.
\end{eqnarray*}
This proves Step 1.

\textbf{Step 2: }\emph{For each $\delta>0$ there exists
  $\epsilon=\epsilon(\delta)>0$ 
such that if $||\nabla\mathcal{A}^H(v,\eta)||\leq\epsilon$, then
$v(t) \in U_\delta$ for every $t \in [0,1]$.}

Denote by $\Gamma_\delta$ the set of smooth paths
$\gamma \in C^\infty([0,1],U_\delta)$ such that 
$|H(\gamma(0))|=\delta$ and $|H(\gamma(1))|=\delta/2$.
For each $x \in U_\delta$ there is a splitting 
$T_x M=T_x H^{-1}(H(x))\oplus T_x ^\perp H^{-1}(H(x))$. 
We denote by $\pi_x$ the projection to the second factor. 
We introduce the number $\epsilon_0=\epsilon_0(\delta)$ by
$$\epsilon_0=\inf_{\gamma \in \Gamma_\delta}\bigg\{\int_0^1
||\pi_{\gamma(t)}(\dot{\gamma}(t))||dt\bigg\}>0.$$
Now assume that $v \in \mathscr{L}$ has the property that
there exist $t_0,t_1 \in \mathbb{R}/\mathbb{Z}$ such that
$|H(v(t_0))| \geq \delta$ and $|H(v(t_1))| \leq \delta/2$.
We claim that
\begin{equation}\label{pen}
||\nabla \mathcal{A}^H(v,\eta)|| \geq \epsilon_0
\end{equation}
for every $\eta \in \mathbb{R}$.
To see that we estimate
\begin{eqnarray*}
||\nabla \mathcal{A}^H(v,\eta)||
&\geq& \sqrt{\int_0^1||\partial_t v-\eta X_H(v)||^2 dt}\\
&\geq& \sqrt{\int_0^1||\pi_v(\partial_t v-\eta X_H(v))||^2 dt}\\
&=& \sqrt{\int_0^1||\pi_v(\partial_t v)||^2 dt}\\
&\geq& \int_0^1||\pi_v(\partial_t v)||dt\\
&\geq& \epsilon_0.
\end{eqnarray*}
This proves (\ref{pen}).
\\
Now assume that $v \in \mathscr{L}$ has the property that
$v(t) \in M \setminus U_{\delta/2}$ for every $t \in [0,1]$. 
In this case we estimate
\begin{equation}\label{npen}
||\nabla \mathcal{A}^H(v,\eta)||
\geq \Bigg|\int_0^1 H(v(t))dt \Bigg| \geq \frac{\delta}{2}
\end{equation}
for every $\eta \in \mathbb{R}$. 
From (\ref{pen}) and (\ref{npen}) Step 2 follows with
$\epsilon<\min\{\epsilon_0, \delta/2\}$.

\textbf{Step 3: }\emph{We prove the proposition.}

Combining Step 1 and Step 2 the proposition follows with
$c_M=2M+\epsilon c_\delta$. \hfill $\square$
\\ \\
Proposition~\ref{mainest} allows us to control the
size of the Lagrange multiplier $\eta$. Our first corollary
considers the case of gradient flow lines. 

\begin{cor}\label{notime}
Assume that $(v,\eta) \in C^\infty(\mathbb{R} \times S^1,V)\times
C^\infty(\mathbb{R},\mathbb{R})$
is a gradient flow line of $\nabla \mathcal{A}^H$ which satisfies
$\lim_{s \to \pm \infty}(v,\eta)(s,\cdot)=(v^\pm,\eta^\pm)(\cdot)
\in \mathrm{crit}(\mathcal{A}^H)$,
where the limit is uniform in the $t$-variable. Then the
$L^\infty$-norm of $\eta$ is bounded uniformly in terms of
a constant $c<\infty$ which only depends on
$\mathcal{A}^H(v^-,\eta^-)$ and $\mathcal{A}^H(v^+,\eta^+)$.
\end{cor}

To prove invariance of our Floer homology under homotopies we
also have to consider the case of $s$-dependent action functionals.
Let $H^-, H^+ \in C^\infty(V)$ be defining Hamiltonians for two  
exact convex hypersurfaces. Consider the
the smooth family of $s$-dependent Hamiltonians $H_s$ defined
as
$$H_s=\beta(s)H^++(1-\beta(s))H^-$$
where $\beta \in C^\infty(\mathbb{R},[0,1])$ is a smooth 
monotone increasing cutoff function such that
$\beta(s)=1$ for $s \geq 1$ and $\beta(s)=0$ for $s \leq 0$. 

\begin{cor}\label{time}
If $\max_{x \in V}|H^+(x)-H^-(x)|$ is small enough, then for each
gradient flow line $(v,\eta)$ of the $s$-dependent action functional
$\mathcal{A}^{H_s}$ which converges at both ends the Lagrange multiplier
$\eta$ is uniformly bounded in terms of the action values at
the end points.
\end{cor}

\textbf{Proof of Corollary~\ref{notime}: } 
Let $\epsilon$ be as in Proposition~\ref{mainest}. 
For $\sigma \in \mathbb{R}$ let $\tau(\sigma) \geq 0$ be defined by
$$\tau(\sigma):=\inf\big\{ \tau \geq 0: ||\nabla\mathcal{A}^H\big
((v,\eta)(\sigma+\tau(\sigma))\big)||
< \epsilon\big\}.$$
We abbreviate the energy of the flow line $(v,\eta)$ by
$$
   E:=\mathcal{A}^H(v^+,\eta^+)-\mathcal{A}^H
(v^-,\eta^-).
$$
We claim that
\begin{equation}\label{away}
\tau(\sigma) \leq \frac{E}{\epsilon^2}.
\end{equation}
To see this we estimate
\begin{eqnarray*}
E&=&\mathcal{A}^H(v^+,\eta^+)-\mathcal{A}^H(v^-,\eta^-)\\
&=&\int_{-\infty}^\infty\frac{d}{ds}\mathcal{A}^H(v,\eta)ds\\
&=&\int_{-\infty}^\infty d\mathcal{A}^H(v,\eta)
\partial_s(v,\eta)ds\\
&=&\int_{-\infty}^\infty \langle \nabla \mathcal{A}^H(v,\eta),
\partial_s(v,\eta)\rangle ds\\
&=&\int_{-\infty}^\infty||\nabla \mathcal{A}^H(v,\eta)||^2 ds\\
&\geq&\int_{\sigma}^{\sigma+\tau(\sigma)}||\nabla \mathcal{A}^H(v,\eta)||^2 ds\\
&\geq&\tau(\sigma)\epsilon^2.
\end{eqnarray*}
This implies (\ref{away}). 
\\
We set
$$M:=\max\{|\mathcal{A}^H(v^+,\eta^+)|,
|\mathcal{A}^H(v^-,\eta^-)|\}.$$ 
Note that since the action increases
along a gradient flow line we have
$$\big|\mathcal{A}^H\big((v,\eta)(\sigma+\tau(\sigma))\big)\big| \leq
M\text{ for all }\sigma\in\R.$$
We deduce from Proposition~\ref{mainest} and the definition of
$\tau(\sigma)$ that
\begin{equation}\label{eta}
\big|\eta\big(\sigma+\tau(\sigma)\big)\big| \leq c_M.
\end{equation}
We set
\begin{equation}\label{cH}
c_H:=\max_{x \in V}|H(x)|.
\end{equation}
We estimate using (\ref{away}), (\ref{eta}), and (\ref{cH})
\begin{eqnarray*}
|\eta(\sigma)|&\leq&\big|\eta\big(\sigma+\tau(\sigma)\big)\big|+
\int_{\sigma}^{\sigma+\tau(\sigma)}|\partial_s \eta(s)|ds\\
&=&\big|\eta\big(\sigma+\tau(\sigma)\big)\big|+
\int_{\sigma}^{\sigma+\tau(\sigma)}\bigg|\int_0^1H(v(s,t)dt\bigg|ds\\
&\leq&c_M+c_H \tau(\sigma)\\
&\leq&c_M+
\frac{c_H E}
{\epsilon^2}.
\end{eqnarray*}
The right hand side is independent of $\sigma$ and hence we get
\begin{equation}\label{bound}
||\eta||_\infty \leq c_M+\frac{c_H E}{\epsilon^2}.
\end{equation}
This proves the corollary. \hfill $\square$
\\ \\
\textbf{Proof of Corollary~\ref{time}: }In the $s$-dependent case we
define the energy as
$$
   E:=\mathcal{A}^{H^+}(v^+,\eta^+)-\mathcal{A}^{H^-}(v^-,\eta^-)
   -\int_{-\infty}^\infty \big(\partial_s \mathcal{A}^{H_s}\big)(v,\eta)
   ds,
$$
where
$$
   \big(\partial_s \mathcal{A}^{H_s}\big)(v,\eta) =
   -\eta\int_0^1\frac{\p H_s}{\p s}(v)dt = -\eta\int_0^1\beta'(s)
   (H^+-H^-)(v)dt.
$$
If we set $c_H:=\max\{c_{H^+},c_{H^-}\}$ and
$\epsilon:=\min\{\epsilon(H^+),\epsilon(H^-)\}$ then (\ref{bound})
can be deduced as in the time-independent case. 
However, $E$ is a priori not bounded anymore because of the
term containing the $s$-derivatives of the action functional.
We  use the abbreviations
$$\Delta:=\mathcal{A}^{H^+}(v^+,\eta^+)-\mathcal{A}^{H^-}(v^-,\eta^-)$$
and 
$$\delta:=\max_{x \in V}|H^+(x)-H^-(x)|$$
and estimate
\begin{eqnarray*}
E&=&\Delta -\int_{-\infty}^\infty (\partial_s \mathcal{A}^{H_s})(v,\eta) ds\\
&=&\Delta +\int_{-\infty}^\infty \beta '(s)\eta(s)\Bigg(\int_0^1
(H^+-H^-)\bigl(v(s,t)\bigr)dt\Bigg)ds\\
&\leq&\Delta +\delta||\eta||_\infty.
\end{eqnarray*}
If we set this estimate into (\ref{bound}) we obtain
$$||\eta||_\infty \leq c_M+\frac{c_H \Delta }{\epsilon^2}+
\frac{c_H \delta}{\epsilon^2}||\eta||_\infty.$$
Now if 
$$\delta<\frac{\epsilon^2}{c_H}$$
we obtain the following uniform $L^\infty$-bound for $\eta$
$$||\eta||_\infty \leq \frac{\epsilon^2 c_M+c_H \Delta }
{\epsilon^2-c_H \delta}.$$
This proves Corollary~\ref{time}. \hfill $\square$
\\ \\
\textbf{Proof of Theorem~\ref{well}: }
As we pointed out at the beginning of this section, we may assume without
loss of generality that $\mathcal{A}^H$ is Morse-Bott. Choose
an additional Morse function $h$ on $\mathrm{crit}(\mathcal{A}^H)$. 
The Floer chain complex is defined in the following
way. $CF(\mathcal{A}^H,h)$ is the
$\mathbb{Z}_2$-vector space consisting of formal sums
$$\xi=\sum_{c \in \mathrm{crit}(h)} \xi_c c$$
where the coefficients $\xi_c \in \mathbb{Z}_2$ satisfy the finiteness
condition
\begin{equation}\label{novikov}
\#\{c \in \mathrm{crit}(h): \xi_c \neq 0,\,\, 
\mathcal{A}^H(c) \leq \kappa\}<\infty
\end{equation}
for every $\kappa \in \mathbb{R}$.
To define the boundary operator, we require 
some compatibility condition of the family of $\omega$-compatible
almost complex structures $J_t$ with the convex structure of $V$ 
at infinity in order
to make sure that our cascades remain in a compact subset of $V$. 
As we remarked in the introduction, completeness implies that there
exists a contact manifold $(M,\alpha_M)$
\footnote{Be careful! Do not confuse the contact manifolds
$(M,\alpha_M)$ and $(\Sigma, \alpha)$.}
such that a neighbourhood of infinity of 
the symplectic manifold  $(V,\omega)$ can be symplectically identified
with 
$(M \times \mathbb{R}_+, d(e^r \alpha_M))$, where
$r$ refers to the coordinate on 
$\mathbb{R}_+=\{r \in \mathbb{R}: r \geq 0\}$. We may assume
without loss of generality that $H$ is constant on
$M \times \mathbb{R}_+$. 
We require the
following conditions on $J_t$ for every $t \in [0,1]$.
\begin{itemize}
 \item For each $x \in M$ we have 
  $J_t(x)\frac{\partial}{\partial r}=R_M$, where $R_M$ 
  is the Reeb vector field on $(M,\alpha_M)$.
 \item $J_t$ leaves the kernel of $\alpha_M$ invariant for
  every $x \in M$.
 \item $J_t$ is invariant under the local 
  half flow $(x,0) \mapsto (x,r)$
  for $(x,r) \in M \times \mathbb{R}_+$.
\end{itemize}
We choose further an additional Riemannian metric
$g_c$ on the critical manifold $\mathrm{crit}(\mathcal{A}^H)$. For
two critical points $c_-, c_+ \in \mathrm{crit}(h)$ we consider
the moduli space of {\em gradient flow lines with cascades}
$\mathcal{M}_{c_-,c_+}(\mathcal{A}^H,h,J,g_c)$ as defined in
Appendix~\ref{app:casc}. . 
For generic choice of $J$ and $g_c$ this moduli space of is a smooth
manifold. We claim that its zero dimensional component
$\mathcal{M}^0_{c_-,c_+}(\mathcal{A}^H,h,J,g_c)$ is
actually compact and hence a finite set. To see that we have to prove
that cascades are compact modulo breaking. Since the support of
$X_H$ lies outside of $M \times \mathbb{R}_+$, the first component of
a gradient flow line which enters $M \times \mathbb{R}_+$ 
will just satisfy the pseudo-holomorphic
curve equation by (\ref{flowline}).
By our choice of the family of almost complex structures the 
convexity condition guarantees that it cannot touch any level set $M
\times \{r\}$ from inside (see~\cite{mcduff}), and since its asymptotics
lie outside of $M \times \mathbb{R}_+$ it has to
remain in the compact set $V \setminus M \times \mathbb{R}_+$ 
all the time. This gives us
a uniform $L^\infty$-bound on the first component. 
Corollary~\ref{notime} implies that the second component remains
bounded, too. Since the symplectic form $\omega$ is exact there are
no nonconstant $J$-holomorphic spheres. This excludes
bubbling and hence the derivatives
of (\ref{flowline}) can be controlled, see \cite{mcduff-salamon}. 
This proves the claim.
\\ 
We now set
$$n(c_-,c_+)=\#\mathcal{M}^0_{c_-,c_+}(\mathcal{A}^H,h,J,g_c) \,\,
\mathrm{mod}\,\,2\in 
\mathbb{Z}_2$$ 
and define the Floer boundary operator
$$\partial \colon CF(\mathcal{A}^H,h) \to CF(\mathcal{A}^H,h)$$
as the linear extension of
$$\partial c=\sum_{c' \in \mathrm{crit}(h)}n(c,c')c'$$
for $c \in \mathrm{crit}(h)$. Again using the fact that
the moduli space of cascades are compact modulo breaking, a well-known
argument in Floer theory shows that 
$$\partial^2=0.$$
We define our Floer homology as usual by
$$HF(\mathcal{A}^H,h,J,g_c)=\frac{\mathrm{ker}\partial}
{\mathrm{im}\partial}.$$
Standard arguments show that $HF(\mathcal{A}^H,h,J,g_c)$ is independent
of the choices of $h$, $J$, and $g_c$ up to canonical isomorphism and
hence $HF(\mathcal{A}^H)$ is well defined. To prove that 
it is invariant under homotopies of $H$ we use Corollary~\ref{notime}
to show that the Floer homotopies which are defined by counting solutions
of the $s$-dependent gradient equation are well defined. This finishes
the proof of Theorem~\ref{well}. \hfill $\square$ 
\\ \\
\textbf{Proof of Theorem~\ref{zero}: }
We consider the following perturbation of $\mathcal{A}^H$. Let
$F\in C^\infty(\mathbb{R}/\mathbb{Z}\times V)$ be a smooth map
such that $F|_{(0,1)\times V}$ has compact support. 
We use the notation $F_t=F(t, \cdot)$ for $t \in \mathbb{R}/\mathbb{Z}$.
Denote by $\phi^t_H$ and $\phi_F^t$ the flows of the Hamiltonian vector
fields of $H$ and $F_t$, respectively. We define
$$\mathcal{A}^H_F \colon \mathscr{L} \times \mathbb{R} \to \mathbb{R}$$
by
$$\mathcal{A}^H_F(v,\eta):=\mathcal{A}^H(v,\eta)-
\int_0^1 F_t \big(\phi^{-t \eta}_H(v(t))\big)dt.$$ 
We further abbreviate
$$\mathfrak{S}(X_H):=\mathrm{cl}\{x \in M: X_H(x) \neq 0\}$$
the support of the Hamiltonian vector field of $H$. 
Theorem~\ref{zero} follows from the following Proposition and a
standard Floer homotopy argument. \hfill $\square$

\begin{prop}\label{nocrit}
Assume that $\phi^1_F(\mathfrak{S}(X_H))\cap \mathfrak{S}(X_H)=\emptyset$.
Then there exists $\tilde{F} \in C^\infty(\mathbb{R}/\mathbb{Z}\times V)$
such that $\tilde{F}|_{(0,1)\times V}$ has compact support and there 
are no critical points of $\mathcal{A}^H_{\tilde{F}}$.
\end{prop}

\textbf{Proof: }
Critical points of $\mathcal{A}^H_F$ are solutions
of the problem
$$
\left. \begin{array}{c}
\partial_t v(t)= \Big(\phi_H^{-\eta t}\Big)^*
X_{F_t}(v(t))+\eta X_H(v(t)),
\quad  t \in \mathbb{R}/\mathbb{Z}, \\
\int_0^1 \big(t\{H,F_t\}(\phi^{-t\eta}_H(v(t))+H(v(t))\big)dt=0
\end{array}
\right\}
$$
with the Poisson bracket given by 
$\{F,H\}=dF(X_H)$. 
Define $w \in C^\infty([0,1],M)$ by
$$w(t):=\phi_H^{-t \eta}(v(t)), \quad t \in [0,1].$$
Then $w$ satisfies
\begin{equation}\label{w}
\left. \begin{array}{c}
\partial_t w(t)=X_{F_t}(w(t)), \quad t \in [0,1],\\
w(1)=\phi^{-\eta}_H(w(0)), \\
\int_0^1 \big(t\{H,F_t\}(w(t))+H(w(t))\big)dt=0.  
\end{array}
\right\}
\end{equation}
For a smooth map $\rho \in C^\infty([0,1],[0,1])$ 
satisfying $\rho(0)=0$ and 
$t \in \mathbb{R}$ set
$$F^\rho_t:=\dot{\rho}(t)F_{\rho(t)} \in C^\infty(V).$$
Note that
$$\phi^t_{F^\rho}=\phi^{\rho(t)}_F.$$
Equations (\ref{w}) for $F$ replaced by $F^\rho$ become
\begin{equation}\label{wnew}
\left. \begin{array}{c}
w(t)=\phi_{F}^{\rho(t)}(w(0)), \\
w(1)=\phi^{-\eta}_H(w(0)), \\
\int_0^1 \big(t\dot{\rho}(t)\{H,F_{\rho(t)}\}(w(t))+H(w(t))\big)dt=0.  
\end{array}
\right\}
\end{equation}
Since the Hamiltonian vector field of $H$ has compact support, there
exists a constant $c$ such that
$$\max_{\substack{x \in V,\\ t \in \mathbb{R}/\mathbb{Z}}}|\{H,F_t\}(x)|\leq c,
\quad \max_{x \in V}|H(x)| \leq c.$$
Using again that the support of the Hamiltonian
vector field is compact together with the fact that
$0$ is a regular value of $H$ we conclude that there exists $\delta>0$
such that
$$\min_{x \in V \setminus \mathfrak{S}(X_H)}|H(x)| = \delta.$$
Choose an $\epsilon>0$ such that
$$\epsilon<\frac{\delta}{2c+\delta}$$
and a smooth function $\rho_\epsilon \in C^\infty([0,1],[0,1])$ such that
$$
\left. \begin{array}{c}
\rho_\epsilon(0)=0\\
\rho_\epsilon(t)=1, \quad t \geq \epsilon, \\
\dot{\rho}_\epsilon(t) \geq 0, \quad t \in [0,1].\\  
\end{array}
\right\}
$$
Proposition~\ref{nocrit} follows now with 
$\tilde{F}=F^{\rho_\epsilon}$ in view of the lemma below.
\hfill $\square$

\begin{lemma}
Assume that $\phi^1_F(\mathfrak{S}(X_H))\cap \mathfrak{S}(X_H)=\emptyset$.
Then there are no solution of (\ref{wnew}) for $\rho=\rho_\epsilon$.
\end{lemma}
\textbf{Proof: }
Let $w$ be a solution of (\ref{wnew}). We first claim that
\begin{equation}\label{claimw}
w(0) \notin \mathfrak{S}(X_H).
\end{equation}
We argue by contradiction and assume that 
$w(0) \in \mathfrak{S}(X_H)$. It follows from the first
equation in (\ref{wnew}) and the assumption of the lemma that
$$w(1)=\phi^{\rho(1)}_F(w(0))=\phi^1_F(w(0))
 \notin \mathfrak{S}(X_H).$$
The definition of $\mathfrak{S}(X_H)$ implies that
$$\phi^{\eta}_H(w(1))=w(1).$$
Combining the above two equations together with the second equation
in (\ref{wnew}) we conclude
$$w(0)=\phi^{\eta}(w(1))=w(1) \notin \mathfrak{S}(X_H).$$
This contradicts the assumption that $w(0) \in \mathfrak{S}(X_H)$
and proves (\ref{claimw}). 
\\
Combining (\ref{claimw}) with the second equation in (\ref{wnew})
we obtain
\begin{equation}\label{w2}
w(1)=w(0) \notin \mathfrak{S}(X_H).
\end{equation}
Using the definition of $\rho_\epsilon$, the first equation 
in (\ref{wnew}), and (\ref{w2}) we get
\begin{equation}\label{w3}
w(t)=\phi^{\rho_\epsilon(t)}_F (w(0))=
\phi^1_F(w(0))=w(1) \notin \mathfrak{S}(X_H), \quad t \geq \epsilon.
\end{equation}
Using the definition of $\mathfrak{S}(X_H)$ and of
$\delta$ we deduce that
\begin{equation}\label{w4}
|H(w(t))|\geq \delta, \quad \{H,F_{\rho_\epsilon(t)}\}
(w(t))=0, \quad t\geq \epsilon.
\end{equation}
Using (\ref{w4}), the definition of $c$ and $\epsilon$, and
the properties of $\rho_\epsilon$  we
estimate
\begin{align*}
&\Bigg|\int_0^1 \big(t\dot{\rho}_\epsilon(t)\{H,F_{\rho_\epsilon(t)}\}
(w(t))+H(w(t))\big)dt\Bigg|\cr
&\geq
-\Bigg|\int_0^\epsilon \big(t\dot{\rho}_\epsilon(t)\{H,F_{\rho_\epsilon(t)}\}
(w(t))+H(w(t))\big)dt\Bigg|\cr
& \ \ 
+\Bigg|\int_\epsilon^1 \big(t\dot{\rho}_\epsilon(t)\{H,F_{\rho_\epsilon(t)}\}
(w(t))+H(w(t))\big)dt\Bigg|
\cr
&\geq
-\Bigg|\int_0^\epsilon (\epsilon c \dot{\rho}_\epsilon(t)+c)dt\Bigg|+
\delta(1-\epsilon)\cr
&=-2c\epsilon+\delta(1-\epsilon)\cr
&>0.
\end{align*}
This contradicts the third equation in (\ref{wnew}). Hence there
are no solutions of (\ref{wnew}),
which proves the lemma. \hfill $\square$

\section{Index computations}\label{sec:index}

In this section we prove Theorem~\ref{compute}. The proof comes down
to the computation of the indices of generators of the Floer chain
complex in the case that $\Sigma$ is the unit cotangent bundle of the
sphere. 

We first have to study the question under which conditions
$HF(\Sigma,V)$ has a $\mathbb{Z}$-grading. Throughout this
section, we make the following assumptions:
\begin{itemize}
\item[(A)]Closed Reeb orbits on $(\Sigma,\alpha)$ are of Morse-Bott
  type (see Section~\ref{sec:floer}). 
\item[(B)]$\Sigma$ is simply connected and $V$ satisfies $I_{c_1}=0$.  
\end{itemize}
Under these assumptions the (transversal) Conley Zehnder index of a
Reeb orbit $v \in C^\infty(S^1,\Sigma)$ can be defined in the
following way. Since $\Sigma$ is simply connected, we can find a 
map $\bar{v} \in C^\infty(D, \Sigma)$ on 
the unit disk $D=\{z \in \mathbb{C}: |z| \leq 1\}$ such that
$\bar{v}(e^{2\pi i t})=v(t)$. Choose a (homotopically unique) 
symplectic trivialization of the symplectic vector bundle
$(\bar{v}^*\xi,\bar{v}^*d\alpha)$. The linearized flow of the
Reeb vector field along $v$ defines a path in the
group $Sp(2n-2,\mathbb{R})$ of symplectic matrices. The Maslov index
of this path \cite{robbin-salamon0} is the {\em (transversal)
  Conley-Zehnder index} $\mu_{CZ}\in\frac{1}{2}\Z$. It is independent
of the choice of the disk $\bar{v}$ due to the assumption $I_{c_1}=0$
on $V$. 

Let $\mathcal{M}$ be the moduli space of all finite energy gradient
flow lines of the action functional
$\mathcal{A}^H$. Since $\mathcal{A}^H$ is Morse-Bott
every finite energy gradient flow line 
$(v,\eta) \in C^\infty(\mathbb{R}\times S^1,V) \times 
C^\infty(\mathbb{R},\mathbb{R})$ converges exponentially
at both ends to critical points 
$(v^\pm,\eta^\pm) \in \mathrm{crit}(\mathcal{A}^H)$ as the flow parameter
goes to $\pm \infty$. The linearization of the gradient flow
equation along any path $(v,\eta)$ in $\mathscr{L}\times \mathbb{R}$ which
converges exponentially to the critical points of $\mathcal{A}^H$
gives rise to an operator $D_{(v,\eta)}^{\mathcal{A}^H}$.
For suitable weighted Sobolev spaces (the weights are needed because
we are in a Morse-Bott situation) the operator
$D_{(v,\eta)}^{\mathcal{A}^H}$ is a Fredholm operator. 
Let $C^-,C^+ \subset \mathrm{crit}(\mathcal{A}^H)$ be
the connected components of the critical manifold of
$\mathcal{A}^H$ containing $(v^-,\eta^-)$ or $(v^+,\eta^+)$
respectively. 
The local virtual dimension of $\mathcal{M}$
at a finite energy gradient flow line is defined to be
\begin{equation}\label{i1}
\mathrm{virdim}_{(v,\eta)}\mathcal{M}:=
\mathrm{ind}D^{\mathcal{A}^H}_{(v,\eta)}+\mathrm{dim}C^-+\mathrm{dim}C^+
\end{equation}
where $\mathrm{ind}D^{\mathcal{A}^H}_{(v,\eta)}$ is the Fredholm index
of the Fredholm operator $D_{(v,\eta)}^{\mathcal{A}^H}$. For
generic compatible almost complex structures, the moduli space of
finite energy gradient flow lines 
is a manifold and the local virtual dimension of the moduli
space at a gradient flow line $(v,\eta)$ corresponds to the
dimension of the connected component of $\mathcal{M}$ containing
$(v,\eta)$. Our first goal is to prove the following index formula. 

\begin{prop}\label{index}
Assume that hypotheses (A) and (B) hold. 
Let $C^-,C^+ \subset \mathrm{crit}(\mathcal{A}^H)$ be
two connected components of the critical manifold of $\mathcal{A}^H$.
Let $(v,\eta) \in C^\infty(\mathbb{R}\times S^1,V) \times
C^\infty(\mathbb{R},\mathbb{R})$ be a gradient flow line
of $\mathcal{A}^H$ which converges at both ends
$\lim_{s \to \pm \infty}(v,\eta)(s) \to (v^\pm,\eta^\pm)$
to critical points of $\mathcal{A}^H$ satisfying
$(v^\pm, \eta^\pm) \in C^\pm$. Choose maps
$\bar{v}^\pm \in C^\infty(D,\Sigma)$ satisfying 
$\bar{v}^\pm(e^{2 \pi i t})=v^\pm(t)$. Then the local virtual
dimension of the moduli space $\mathcal{M}$
of finite energy gradient flow lines
of $\mathcal{A}^H$ at $(v,\eta)$ is given by
\begin{equation}\label{virdim}
\mathrm{virdim}_{(v,\eta)}\mathcal{M}=
\mu_{CZ}(v^+)-\mu_{CZ}(v^-)
+2c_1(\bar{v}^-\#v\#\bar{v}^+)
+\frac{\mathrm{dim}C^-+\mathrm{dim}C^+}{2}
\end{equation}
where $\bar{v}^-\#v\#\bar{v}^+$ is the sphere obtained by
capping the cylinder $v$ with the disks $\bar{v}^+$ and
$\bar{v}^-$, and $c_1=c_1(TV)$. 
\end{prop}

The proof is based on a discussion of

{\bf Spectral flows. }It is shown in \cite{robbin-salamon} that the
Fredholm index of $D^{\mathcal{A}^H}_{(v,\eta)}$ can be computed via
the {\em spectral flow} $\mu_\spec$ (see Appendix~\ref{sec:spec}) of
the Hessian $\Hess_{\mathcal{A}^H}$ along $(v,\eta)$ by the formula
\begin{equation}\label{i2}
\mathrm{ind}D^{\mathcal{A}^H}_{(v,\eta)} =
\mu_\spec\Bigl(\Hess_{\mathcal{A}^H}(v,\eta)\Bigr).  
\end{equation}
Our proof compares the
spectral flow of the Hessian of $\mathcal{A}^H$
with the spectral flow of the action functional of classical mechanics
which can be computed via the Conley-Zehnder indices. 
For a {\em fixed} Lagrange multiplier $\eta \in \mathbb{R}$
the action functional of classical mechanics arises as
$$
   \mathcal{A}^H_\eta := \mathcal{A}^H(\cdot,\eta) \colon \mathscr{L}
\to \mathbb{R}.
$$ 

Assume first that the periods $\eta^\pm$ of the Reeb orbits $v^\pm$
are nonzero. We begin by homotoping the action functional
$\mathcal{A}^H$ via Morse-Bott functionals with fixed critical manifold
to an action functional $\mathcal{A}^{H^1}$ which satisfies
the assumptions of (the infinite dimensional analogue of)
Lemma~\ref{varlag}. There exists a neighbourhood $U \subset V$ of
$\Sigma$ and an $\epsilon>0$ such that $U$ is
symplectomorphic to 
$\big(\Sigma \times (-\epsilon, \epsilon), d (e^r \alpha)\big)$
where $r$ is the coordinate on $(-\epsilon, \epsilon)$. 
Since $\AA^H$ is Morse-Bott and the Hamiltonian vector field $X_H(x)$
for $x \in \Sigma$ equals the Reeb vector field $R(x)$, there exists a 
homotopy $H^s$ for $s \in [0,1]$ which satisfies the following
conditions: 
\begin{itemize}
 \item $H ^0=H$.
 \item $X_{H^s}(x)=R(x)$ for $x \in \Sigma$ and $s \in [0,1]$.
 \item There exist neighbourhoods $U^\pm\subset U$ of the critical
  manifolds $C^\pm$ and functions $h_\pm \in C^\infty((-\epsilon,\epsilon))$
  satisfying $h_\pm(0)=0$, $h_\pm'(0)=1$, 
  $h_\pm''(0) \neq 0$, and $h_\pm'(r) \neq 0$ for
  $r \in (-\epsilon, \epsilon)$ such that
  $H^1(x,r)=h_\pm(r)$ for $(x,r)\in
  U^\pm\subset\Sigma\times(-\epsilon,\epsilon)$. 
 \item $\AA^{H^s}$ is Morse-Bott for all $s\in[0,1]$. 
\end{itemize}
Here the signs of $h_\pm''(0)$ are determined by the second derivatives
of $H$ in the direction transverse to $\Sigma$ along $C^\pm$. 
Since $\mathcal{A}^H$ can be homotoped to $\mathcal{A}^{H^1}$ 
via Morse-Bott action functionals with fixed critical manifold, we
obtain
\begin{equation}\label{i3}
\mu_\spec\Bigl(\Hess_{\mathcal{A}^H}(v,\eta)\Bigr) =
\mu_\spec\Bigl(\Hess_{\mathcal{A}^{H^1}}(v,\eta)\Bigr).
\end{equation}
If $(v_0,\eta_0) \in C^\infty(S^1,\Sigma\cap U^\pm) \times \mathbb{R}$
is a critical point of
$\mathcal{A}^H$, then $(v_0,\eta_0)$ is also a critical point of
$\mathcal{A}^{H^1}$. Moreover, the family
$(v_\rho,\eta_\rho) \in C^\infty(S^1,U)\times \mathbb{R}$ 
given by
$$v_\rho(t)=(v_0(t),h^{-1}(-\rho)), \quad 
\eta_\rho=\frac{\eta_0}{h'(h^{-1}(-\rho))}$$
consists of critical points for the family of action functionals
$\mathcal{A}^{H^1,\rho} \colon \mathscr{L}\times \mathbb{R} \to \mathbb{R}$
given for $(v,\eta) \in \mathscr{L} \times \mathbb{R}$ by
$$\mathcal{A}^{H^1,\rho}(v,\eta) := \int v^*\lambda-
\eta\bigg(\int_0^1H^1(v(t))dt+\rho\bigg).$$
Note that
$$
   \p_\rho\eta_\rho|_{\rho=0} =
   -\frac{\eta_0h_\pm''(0)}{h_\pm'(0)^2}. 
$$
Hence for $\eta_0=\eta^\pm\neq 0$ the hypotheses of Lemma~\ref{varlag} are
satisfied. It follows from Theorem~\ref{lagspec} and
Lemma~\ref{varlag} that the spectral flow can be expressed in terms of
the spectral flow of the action functional of classical mechanics plus
a correction term accounting for the second derivatives of $H$
transversally to $\Sigma$ as
\begin{equation}\label{i4}
\mu_\spec\Bigl(\Hess_{\mathcal{A}^{H^1}}(v,\eta)\Bigr)=
\mu_\spec\Bigl(\Hess_{\mathcal{A}^{H^1}_\eta}(v)\Bigr)+\frac{1}{2}
\bigg( \mathrm{sign}(\eta^-\cdot h_-''(0))-
\mathrm{sign}(\eta^+ \cdot h_+''(0))\bigg).
\end{equation}
It follows from a theorem due to Salamon and Zehnder
\cite{salamon-zehnder} that the spectral flow of the Hessian of
$\mathcal{A}^{H^1}_\eta$ can be computed via Conley-Zehnder
indices. However, the Conley-Zehnder 
indices in the Salamon-Zehnder theorem are not the
(transversal) Conley-Zehnder indices explained above, 
but the Maslov index of the linearized flow of the Reeb vector
field on the whole tangent space of $V$ and not just on the contact
hyperplane. For a Reeb orbit $v$ we will denote this second {\em
  (full) Conley-Zehnder index} by $\hat{\mu}_{CZ}(v)$. Note that
$\hat{\mu}_{CZ}(v)$ depends 
on the second derivatives of $H$ transversally to $\Sigma$ while
$\mu_{CZ}(v)$ does not. Another complication is that we are in a Morse-Bott
situation and we have to adapt the Salamon-Zehnder theorem to this situation.
Formula (\ref{extsf}) defines the spectral flow also for Morse-Bott
situations. To adopt the Conley-Zehnder indices to the Morse-Bott situation
observe that in a symplectic trivialization the linearized flow
of the Reeb vector field can be expressed as a solution of an ordinary
differential equation
$$\dot{\Psi}(t)=J_0 S(t)\Psi(t), \quad \Psi(0)=\mathrm{id},$$
where $t \mapsto S(t)=S(t)^T$ is a smooth path of symmetric matrices. 
For a real number $\delta$ we define $\Psi_\delta$ as the
solution of 
$$\dot{\Psi}_\delta(t)=J_0 \big(S(t)-\delta\cdot \mathrm{id}\big)\Psi_\delta(t), 
\quad \Psi_\delta(0)=\mathrm{id},$$
and set $\mu^{\delta}_{CZ}(v)$, respectively $\hat{\mu}^{\delta}_{CZ}(v)$
as the Conley-Zehnder index of $\Psi_{\delta}$ where in the first case
we restrict $\Psi_\delta$ to the contact hyperplane and in the second
case we consider it on the whole tangent space. We put
$$\mu^+_{CZ}(v):=\lim_{\delta \searrow 0}\mu^{\delta}_{CZ}(v), \quad
\mu^-_{CZ}(v):=\lim_{\delta \searrow 0}\mu^{-\delta}_{CZ}(v)$$
and analoguously $\hat{\mu}^+_{CZ}(v)$ and $\hat{\mu}^-_{CZ}(v)$. Note
that while $\hat{\mu}_{CZ}(v)$ and $\mu_{CZ}(v)$ are half-integers,
$\hat{\mu}^\om_{CZ}(v)$ and $\mu^\pm_{CZ}(v)$ are actually integers.  
We are now in position to state the theorem of Salamon and Zehnder.
\begin{thm}[Salamon-Zehnder~\cite{salamon-zehnder}]\label{diedi}
The spectral flow of the Hessian of $\mathcal{A}^{H^1}_\eta$ is
given by
$$\mu(H_{\mathcal{A}^{H^1}_\eta}(v))=\hat{\mu}^+_{CZ}(v^+)-
\hat{\mu}^-_{CZ}(v^-)+2c_1(\bar{v}^-\#v\#\bar{v}^+).$$
\end{thm}

{\bf Relations between Conley-Zehnder indices. }
The following two lemmata relate the different Conley-Zehnder indices to
each other. 
\begin{lemma}\label{coz1}
For a Reeb orbit $v$ with period $\eta\neq 0$, viewed as a 1-periodic
orbit of the Hamiltonian vector field of $\eta H$, we have  
$$\hat{\mu}^\pm_{CZ}(v)=\mu^\pm_{CZ}(v)+\frac{1}{2}\bigg(
\mathrm{sign}\big(\eta h''(0)\big)\mp 1\bigg).$$
\end{lemma}
\textbf{Proof: } By the product property \cite{salamon} of the
Conley-Zehnder index the difference of
$\hat{\mu}^\pm_{CZ}(v)$ and $\mu^\pm_{CZ}(v)$ is given by
the Conley Zehnder index of the linearized flow of the Hamiltonian vector
field restricted to the symplectic orthogonal complement $\xi^\om$ of
the contact hyperplane in the tangent space of $V$. With respect
to the trivialization $\mathbb{C} \to \xi^\omega(v(t))$ given
by $x+iy \mapsto (x \cdot \nabla H(v(t))+y \cdot X_H(v(t)))$ for $t \in
S^1$, 
the linearized flow of the Hamiltonian vector field is given by
$$\Psi(t)=\left(\begin{array}{cc}
1 & 0\\
t\eta h''(0) & 1
\end{array}\right).$$
Abbreviate $a := \eta \cdot h''(0)$. 
A computation shows that
$$
   \Psi_\delta(t) = e^{\delta(a-\delta)t^2}
   \left(\begin{array}{cc}
      1 & \delta t \\
      (a-\delta) t & 1
   \end{array}\right). 
$$ 
Recall \cite{salamon} that the Conley-Zehnder index can be computed
in terms of crossing numbers, where a number $t \in [0,1]$ is called
a crossing if $\mathrm{det}(\mathrm{id}-\Psi_{\delta}(t))=0$. 
The formula above shows that for $\delta$ small enough the only
crossing happens at zero. Hence by~\cite{salamon} the
Conley-Zehnder index is given by
$$\mu_{CZ}(\Psi_\delta)=\frac{1}{2}\mathrm{sign}
\left(\begin{array}{cc}
a-\delta & 0\\
0 &-\delta
\end{array}\right).$$
If $|\delta|<|a|$ we obtain
$$\mu_{CZ}(\Psi_\delta)=\frac{1}{2}\bigg(
\mathrm{sign}(a)-\mathrm{sign}(\delta)\bigg)=
\frac{1}{2}\bigg(\mathrm{sign}\big(\eta 
h''(0)\big)-\mathrm{sign}(\delta)\bigg)$$
and hence
$$\hat{\mu}_{CZ}^\pm(v)-\mu_{CZ}^\pm(v)=
\frac{1}{2}\bigg(\mathrm{sign}\big(\eta h''(0)\big)\mp 1\bigg).$$
This proves the lemma. \hfill $\square$
\begin{lemma}\label{coz2}
Let $v$ be a Reeb orbit with period $\eta\neq 0$ and $C_v$ the
component of the critical manifold of $\mathcal{A}^{H^1}$ which
contains $v$. Then
$$\hat{\mu}_{CZ}(v)=\hat{\mu}^\pm_{CZ}(v)\pm \frac{\mathrm{dim}C_v}{2},
\quad \mu_{CZ}(v)=\mu_{CZ}^\pm(v)\pm \frac{\mathrm{dim}C_v-1}{2}.$$
\end{lemma}
\textbf{Proof: }Obviously
\begin{equation}\label{cz1}
\hat{\mu}^-_{CZ}(v)-\hat{\mu}^+_{CZ}(v)=\mathrm{dim}C_v,
\quad 
\mu^-_{CZ}(v)-\mu^+_{CZ}(v)=\mathrm{dim}C_v-1.
\end{equation}
The reason for the minus one in the second formula is that
the transversal Conley-Zehnder index only takes into account
the critical manifold of $\mathcal{A}^{H^1}$ modulo the
$S^1$-action given by the Reeb vector field. The Conley-Zehnder index
can be interpreted as intersection number of a path of Lagrangian
subspaces with the Maslov cycle, see \cite{robbin-salamon0}.
Under a small perturbation the intersection number can only
change at the initial and endpoint. Since the Lagrangian subspace
at the initial point is fixed it will change only at the endpoint.
There the contribution is given by half of the crossing number
which equals $\mathrm{dim}C_v$ in the case one considers the
Conley-Zehnder index on the whole tangent space respectively
$\mathrm{dim}C_v-1$ if one considers the Conley-Zehnder index only
on the contact hyperplane. In particular,
\begin{equation}\label{cz2}
|\hat{\mu}_{CZ}(v)-\hat{\mu}_{CZ}^\pm(v)| \leq \frac{\mathrm{dim}C_v}{2},
\quad |\mu_{CZ}(v)-\mu^\pm_{CZ}(v)| \leq \frac{\mathrm{dim}C_v-1}{2}.
\end{equation} 
Comparing (\ref{cz1}) and (\ref{cz2}) the lemma follows. \hfill $\square$
\\ \\
\textbf{Proof of Proposition~\ref{index}: }
We first assume that $\eta^-$ and $\eta^+$ are nonzero. 
Combining the theorem of Salamon and Zehnder (Theorem~\ref{diedi}) with
Lemma~\ref{coz1} and Lemma~\ref{coz2}
we obtain
\begin{eqnarray*}
\mu(H_{\mathcal{A}^{H^1}_\eta}(v))&=&
\hat{\mu}^+_{CZ}(v^+)-\hat{\mu}^-_{CZ}(v^-)+
2c_1(\bar{v}^-\#v\#\bar{v}^+)\\
&=&\mu^+_{CZ}(v^+)-\mu^-_{CZ}(v^-)+
2c_1(\bar{v}^-\#v\#\bar{v}^+)-1\\
& &+\frac{1}{2}\bigg(\mathrm{sign}(\eta^+ \cdot h_+''(0))
-\mathrm{sign}(\eta^- \cdot h+-''(0))\bigg)\\
&=&\mu_{CZ}(v^+)-\mu_{CZ}(v^-)+
2c_1(\bar{v}^-\#v\#\bar{v}^+)-\frac{\mathrm{dim}C^-+\mathrm{dim}C^+}{2}\\
& &+\frac{1}{2}\bigg(\mathrm{sign}(\eta^+ \cdot h_+''(0))
-\mathrm{sign}(\eta^- \cdot h_-''(0))\bigg).
\end{eqnarray*}
Combining this equality with (\ref{i1}), (\ref{i2}), (\ref{i3}), and
(\ref{i4}) we compute
\begin{eqnarray*}
\mathrm{virdim}_{(v,\eta)}\mathcal{M}
&=&\mathrm{ind}D^{\mathcal{A}^H}_{(v,\eta)}+\mathrm{dim}C^-+\mathrm{dim}C^+\\
&=&\mu(H_{\mathcal{A}^{H}}(v,\eta))+\mathrm{dim}C^-+\mathrm{dim}C^+\\
&=&\mu(H_{\mathcal{A}^{H^1}}(v,\eta))+\mathrm{dim}C^-+\mathrm{dim}C^+\\
&=&\mu(H_{\mathcal{A}^{H^1}_\eta}(v))
+\frac{1}{2}\bigg(\mathrm{sign}(\eta^-\cdot h_-''(0))-
\mathrm{sign}(\eta^+\cdot h_+''(0))\bigg)
\\
& &+\mathrm{dim}C^-+\mathrm{dim}C^+\\
&=&\mu_{CZ}(v^+)-\mu_{CZ}(v^-)
+2c_1(\bar{v}^-\#v\#\bar{v}^+)\\
& &+\frac{\mathrm{dim}(C^-)+\mathrm{dim}(C^+)}{2}.
\end{eqnarray*}
This proves the proposition for the case where the periods of the
asymptotic Reeb orbits are both nonzero. To treat also the case where
one of the asymptotic Reeb orbits is constant we consider
the following involution on the loop space $\mathscr{L}$
$$I(v)(t)=v(-t), \quad v \in \mathscr{L},\,\, t\in S^1.$$
We extend this involution to an involution on $\mathscr{L}\times \mathbb{R}$
which we denote by abuse of notation also by $I$ and which is given by
$$I(v,\eta)=(I(v),-\eta), \quad (v,\eta) \in \mathscr{L}\times \mathbb{R}.$$
The action functional $\mathcal{A}^H$ transforms under the involution $I$
by
$$\mathcal{A}^H(I(v,\eta))=-\mathcal{A}^H(v,\eta),\quad
(v,\eta) \in \mathscr{L}\times \mathbb{R}.$$
In particular, the restriction of the involution $I$ to the
critical manifold of $\mathcal{A}^H$ induces an involution
on $\mathrm{crit}(\mathcal{A}^H)$ and the fixed points of this involution
are the constant Reeb orbits. 
\\
We consider now a finite energy gradient flow line
$(v,\eta) \in C^\infty(\mathbb{R}\times S^1,V)\times 
C^\infty(\mathbb{R},\mathbb{R})$ of the action functional
$\mathcal{A}^H$ whose right end $(v^+,\eta^+)$ is a constant Reeb orbit
and whose left end $(v^-,\eta^-)$ is a nonconstant Reeb orbit. 
For the path $(v,\eta)$ in $\mathscr{L}\times \mathbb{R}$ we
consider the path $(v,\eta)_I=(v_I,\eta_I)$ 
in $\mathscr{L}\times \mathbb{R}$ defined
by $(v,\eta)_I(s)=I(v,\eta)(-s)$ for $s \in \mathbb{R}$. 
The path $(v,\eta)_I$ goes from $(v^+,\eta^+)$ to 
$I(v^-,\eta^-)$ and gluing the paths $(v,\eta)$ and $(v,\eta)_I$ together
we obtain a path $(v,\eta)\#(v,\eta)_I$ from
$(v^-,\eta^-)$ to $I(v^-,\eta^-)$. The Fredholm indices of the different
paths are related by
$$\mathrm{ind}D^{\mathcal{A}^H}_{(v,\eta)}=
\mathrm{ind}D^{\mathcal{A}^H}_{(v,\eta)_I}, \quad
\mathrm{ind}D^{\mathcal{A}^H}_{(v,\eta)\#(v,\eta)_I}=
\mathrm{ind}D^{\mathcal{A}^H}_{(v,\eta)}+
\mathrm{ind}D^{\mathcal{A}^H}_{(v,\eta)_I}
+\mathrm{dim}C^+.$$
From this we compute, using (\ref{virdim}) for the case of
nonconstant Reeb orbits and the equality
$\mu_{CZ}(I(v^\pm))=-\mu_{CZ}(v^\pm)$,  
\begin{eqnarray*}
\mathrm{ind}D^{\mathcal{A}^H}_{(v,\eta)}&=&
\frac{1}{2}\cdot\mathrm{ind}D^{\mathcal{A}^H}_{(v,\eta)\#(v,\eta)_I}
-\frac{\mathrm{dim}C^+}{2}\\
&=&\frac{1}{2}\bigg(\mu_{CZ}(I(v^-))-\mu_{CZ}(v^-)+
2c_1(\bar{v}^-\#v\#v_I\#I\bar{v}^+)\\
& &-\frac{\mathrm{dim}C^-+\mathrm{dim}IC^-}{2}\bigg)
-\frac{\mathrm{dim}C^+}{2}\\
&=&-\mu_{CZ}(v^-)+2c_1(\bar{v}^-\#v)-
\frac{\mathrm{dim}C^-+\mathrm{dim}C^+}{2}, 
\end{eqnarray*}
from which we deduce (\ref{virdim}) using (\ref{i1}). This proves
the proposition for the case of gradient flow lines whose left end is
a constant Reeb orbit. The case of gradient flow lines whose right end
is constant can be deduced in the same way or by considering the coindex.
This finishes the proof of the Proposition~\ref{index}. \hfill $\square$
\\ \\
In order to define a $\mathbb{Z}$-grading on $HF(\Sigma,V)$
we need that the local virtual dimension just depends on the asymptotics of 
the finite energy gradient flow line. By (\ref{virdim}) this is the
case if $I_{c_1}=0$ on $V$. In this case the local virtual dimension
is given by
\begin{equation}\label{virdim2}
\mathrm{virdim}_{(v,\eta)}\mathcal{M}=
\mu_{CZ}(v^+)-\mu_{CZ}(v^-)
+\frac{\mathrm{dim}C^-+\mathrm{dim}C^+}{2}.
\end{equation}
In order to deal with the third term it is useful 
to introduce the following index for the Morse function 
$h$ on $\mathrm{crit}(\mathcal{A}^H)$. We define the {\em 
  signature index} $\mathrm{ind}^\sigma_h(c)$ of a critical point $c$
of $h$ to be 
$$\mathrm{ind}^\sigma_h(c):=-\frac{1}{2}\mathrm{sign}(\Hess_h(c)),$$
see Appendix~\ref{app:casc}. 
The half signature index is related to the {\em Morse index} 
$\mathrm{ind}^m_h(c)$, given by the
number of negative eigenvalues of $\Hess_h(c)$ counted with
multiplicity, by
\begin{equation}\label{signind}
\mathrm{ind}^\sigma_h(c) = -\mathrm{ind}^m_h(c)-\frac{1}{2}
\mathrm{dim}_c\big(\mathrm{crit}(\mathcal{A}^H)\big).
\end{equation}
We define a {\em grading} $\mu$ on $CF_*(\mathcal{A}^H,h)$ by
$$\mu(c):=\mu_{CZ}(c)+\mathrm{ind}^\sigma_h(c).$$
By considering the case of nondegenerate closed Reeb orbits, one sees
that $\mu$ takes values in the set $\frac{1}{2}+\Z$, so it is indeed a
$\Z$-grading (shifted by $\frac{1}{2}$). 
Using equation (\ref{virdim2}), it is shown in Appendix~\ref{app:casc}
that the Floer boundary operator $\partial$ has degree $-1$ with
respect to this grading. Hence we get a
$\mathbb{Z}$-grading on the homology $HF_*(\Sigma,V)$. 
\\ \\
\textbf{Proof of Theorem~\ref{compute}: }To prove
Theorem~\ref{compute} we use the fact that
the chain groups underlying the Floer homology $HF_*(\mathcal{A}^H)$
only depend on $(\Sigma,\alpha)$ and not on the embedding
of $\Sigma$ into $V$. We show that for the unit cotangent bundle
$S^*S^n$ for $n \geq 4$ the Floer homology equals the chain complex. 
More precisely, we choose the standard round metric on $S^n$
normalized such that all geodesics are closed with
minimal period one. For this choice assumption (A) 
is satisfied. The critical manifold of $\mathcal{A}^H$ consists of
$\mathbb{Z}$ copies 
of $S^*S^n$, where $\mathbb{Z}$ corresponds to the period of the
geodesic. There is a Morse function $h_0$ on $S^*S^n$ with precisely 4
critical points and zero boundary operator (with $\Z_2$-coefficients!)
whose Morse homology satisfies
$$
HM_k(S^*S^n;\mathbb{Z}_2)=CM_k(h_0;\mathbb{Z}_2)=\left\{
\begin{array}{cc}
\mathbb{Z}_2 & k \in \{0,n-1,n,2n-1\}\\
0 & \mathrm{else}.
\end{array}
\right.
$$
Let $h$ be the Morse function on the critical manifold which coincides
with $h_0$ on each connected component. The chain complex is
generated by
$$\mathrm{crit}(h) \cong \mathbb{Z} \times \mathrm{crit}(h_0).$$
A closed geodesic $c$ is also a critical point of the energy functional
on the loop space. The {\em index} $\mathrm{ind}_E(c)$ of a closed geodesic is 
defined to be the Morse index of the energy functional at the geodesic
and the {\em nullity} $\nu(c)$ is defined to be the dimension
of the connected component of the critical manifold of the energy 
functional which contains the geodesic minus one.
The (transverse) Conley-Zehnder index of a closed geodesic is given by 
\begin{equation}\label{morse}
   \mu_{CZ}(c)=\mathrm{ind}_E(c)+\frac{\nu(c)}{2}.
\end{equation}
This is proved in \cite{duistermaat, weber} for nondegenerate
geodesics; the degenerate case follows from the nondegenerate one
using a the averaging property of the Conley-Zehnder index
(Lemma~\ref{coz2}).  
By the Morse index theorem, see \cite{morse} or 
\cite[Theorem 2.5.14]{klingenberg}, the index of a geodesic is
given by the number of conjugate points counted with multiplicity
plus the concavity. The latter one vanishes for the standard round
metric on $S^n$, since each closed geodesic has a variation
of closed geodesics having the same length \cite{ziller}.
\\
Using the Morse index theorem and equations~\eqref{morse}
and~\eqref{signind}, we compute the index of $(m,x) \in
\mathbb{Z}\times \mathrm{crit}(h_0)$:
\begin{eqnarray*}
\mu(m,x)&=&\mu_{CZ}(m,x)+\mathrm{ind}^{\sigma}_h(x)\\
&=&\mathrm{ind}_E(m,x)+\frac{\nu(m,x)}{2}+\mathrm{ind}^\sigma_h(x)\\
&=&(2m-1)(n-1)+\frac{2n-2}{2}+\mathrm{ind}^\sigma_h(x)\\
&=&2m(n-1)+\mathrm{ind}^m_h(x)-\frac{2n-1}{2}.
\end{eqnarray*}
It follows from Lemma~\ref{period} that the action satisfies
$$\mathcal{A}^H(m,x)=m.$$
In order to have a gradient flow line of $\mathcal{A}^H$ from a
critical point $(m_1,x_1)$ to a critical point $(m_2,x_2)$ we need 
$$
   \mathcal{A}^H(m_2,x_2)-\mathcal{A}^H(m_1,x_1)=m_2-m_1>0
$$
and
$$
   \mu(m_2,x_2)-\mu(m_1,x_1) = 2(m_2-m_1)(n-1)+(i_2-i_1) = 1
$$
for $i_1,i_2\in\{0,n-1,n,2n-1\}$, which is impossible if $n \geq 4$. 
Hence there are no gradient flow lines, so the Floer homology
equals the chain complex. This proves Theorem~\ref{compute}. \hfill
$\square$

\appendix

\section{Morse-Bott homology}\label{app:casc}

We briefly indicate in this appendix how to define Morse-Bott homology
by the use of gradient flow lines with cascades. More details of this
approach can be found in \cite[Appendix A]{frauenfelder}. We begin
with the finite dimensional situation. Let $M$ be
a manifold and $f \in C^\infty(M)$ a {\em Morse-Bott function},
i.e.\,the critical set $\mathrm{crit}(f)$ is a manifold and 
$$
   T_x \mathrm{crit}(f)=\mathrm{ker}\Hess_f(x), \quad x \in
   \mathrm{crit}(f),
$$
where $\Hess_f(x)$ denotes the Hessian of $f$ at $x$. We then choose
an additional Morse function $h$ on $\mathrm{crit}(f)$. The chain
group for Morse-Bott homology is the $\mathbb{Z}_2$-vector space given by
$$CM(f,h):=\mathrm{crit}(h)\otimes \mathbb{Z}_2.$$
Morse-Bott homology should also be definable over the integers via
the cascade approach, but this is nowhere written down. The boundary
operator is defined by counting gradient flow lines with cascades 
between two critical points of $h$ which are indicated by the following
picture.

\medskip
\begin{center}\begin{picture}(0,0)%
\includegraphics{f1.pstex}%
\end{picture}%
\setlength{\unitlength}{3947sp}%
\begingroup\makeatletter\ifx\SetFigFont\undefined%
\gdef\SetFigFont#1#2#3#4#5{%
  \reset@font\fontsize{#1}{#2pt}%
  \fontfamily{#3}\fontseries{#4}\fontshape{#5}%
  \selectfont}%
\fi\endgroup%
\begin{picture}(3624,3312)(2839,-4636)
\put(4726,-3211){\makebox(0,0)[lb]{\smash{\SetFigFont{12}{14.4}{\familydefault}{\mddefault}{\updefault}{\color[rgb]{0,0,0}$\dot{y_1}=-\nabla h(y_1)$}%
}}}
\put(3001,-2311){\makebox(0,0)[lb]{\smash{\SetFigFont{12}{14.4}{\rmdefault}{\mddefault}{\updefault}{\color[rgb]{0,0,0}$\dot{x_1}=-\nabla f(x_1)$}%
}}}
\end{picture}
\end{center}
\medskip

A {\em gradient flow line with cascades} starts with a gradient flow line of 
$h$ on $\mathrm{crit}(f)$ which converges at its negative asymptotic end to a 
critical point of $h$. In finite time this gradient flow line
meets the asymptotic end of a gradient flow line of the Morse-Bott function
$f$. We refer to this gradient flow line of $f$ as the first cascade. 
The cascade converges at its positive end again to a point in 
$\mathrm{crit}(f)$. There the flow continuous with the gradient flow of
$h$ on $\mathrm{crit}(f)$. After finite time a second cascade might appear
but having passed through finitely many cascades
we finally end up with a gradient flow line of $h$
which we follow until it converges asymptotically to a critical point of $h$. 
Gradient flow lines with zero cascades are also allowed. They correspond
to ordinary Morse flow lines of the gradient of $h$ on the manifold
$\mathrm{crit}(f)$. For a formal definition of gradient flow lines
with cascades we refer to \cite{frauenfelder}.

We finally discuss the grading for Morse-Bott homology. If $f$ is Morse
then the following two index conventions are often used. Either Morse homology
is graded by the {\em Morse index} $\mathrm{ind}_f^m$, 
given by the dimensions of the negative
eigenspaces, or by the {\em signature index} 
$$
\mathrm{ind}^\sigma_f(x):=-\frac{1}{2}\mathrm{sign}\Hess_f(x),\quad
x \in \mathrm{crit}(f),
$$
where $\mathrm{sign}$ denotes the signature of the quadratic form
given by the difference of the dimensions of the positive and negative
eigenspaces. The two indices are related by the following global shift
\begin{equation}\label{shift}
   \mathrm{ind}^\sigma_f=\mathrm{ind}^m_f-\frac{1}{2}\mathrm{dim}(M).
\end{equation}
In particular, if $M$ is even dimensional then the signature index 
is integer valued and if $M$ is odd dimensional then it is 
half integer valued. The signature index plays an important role
in Floer's semi-infinite dimensional Morse theory. There the 
stable and unstable manifolds are both infinite
dimensional and hence the Morse index is infinite. The grading given
is given by the  Maslov index which can be interpreted as a signature
index as explained in \cite{robbin-salamon1, robbin-salamon2}.  
\\
Both the Morse index and the signature index can be defined in
the same way also for a Morse-Bott function $f$. The corresponding
indices for a pair $(f,h)$ consisting of a Morse-Bott function $f$
and a Morse function $h$ on $\mathrm{crit}(f)$ are defined
by taking sums
$$\mathrm{ind}^m_{f,h}(x):=\mathrm{ind}^m_f(x)+\mathrm{ind}^m_h(x),
\quad \mathrm{ind}^\sigma_{f,h}(x):=\mathrm{ind}^\sigma_f(x)+
\mathrm{ind}^\sigma_h(x),\quad x \in \mathrm{crit}(h).$$
The shift formula (\ref{shift}) continues to hold for 
these indices in Morse-Bott theory. 

Consider now gradient flow lines with $k$ cascades between components
$C_0,\dots,C_k$ of $\mathrm{crit}(f)$, starting at a critical point
$x^+$ of $h$ on $C^+=C_k$ and ending at a critical point $x^-$ of $h$ on
$C^-=C_0$. For generic metric, their moduli space (divided by the
$\R$-actions on the cascades) 
$\MM(x^-,x^+;C_0,\dots,C_k)$ is a manifold of dimension 
\begin{align*}
   \dim\MM(x^-,x^+;C_0,\dots,C_k) &=
   \ind_h^m(x^+)-\ind_h^m(x^-)-1 \cr
   &\ \ +\sum_{i=1}^k\Bigl(\dim\MM(C_{i-1},C_i)
   -\dim C_i\Bigr),  
\end{align*}
where $\MM(C_{i-1},C_i)$ is the moduli space of gradient flow lines of
$f$ from $C_{i-1}$ to $C_i$ ({\em not} divided by the
$\R$-action). From
\begin{equation}\label{dim1}
   \dim\MM(C_{i-1},C_i) = \ind_f^m(C_i)-\ind_f^m(C_{i-1})+\dim C_i
\end{equation}
we obtain the dimension formula in terms of Morse indices
\begin{align}\label{dim-casc1}
   \dim\MM(x^-,x^+;C_0,\dots,C_k) 
   &=\ind_h^m(x^+)-\ind_h^m(x^-)-1 + \ind_f^m(C^+)-\ind_f^m(C^-) \cr 
   &= \ind_{f,h}^m(x^+)-\ind_{f,h}^m(x^-)-1.
\end{align}
On the other hand, in the Morse-Bott case the Morse and signature
indices of a critical component $C$ are related by
$$
   \ind_f^\sigma(C) = \ind_f^m(C)-\frac{1}{2}(\dim M-\dim C). 
$$
Inserting this in equation~\eqref{dim1} yields
\begin{equation}\label{dim2}
   \dim\MM(C_{i-1},C_i) = \ind_f^\sigma(C_i)-\ind_f^\sigma(C_{i-1})+
   \frac{\dim C_i+\dim C_{i-1}}{2},
\end{equation}
which in turn yields the dimension formula in terms of signature indices
\begin{align}\label{dim-casc2}
   \dim\MM(x^-,x^+;C_0,\dots,C_k) 
   &= \ind_h^m(x^+)-\ind_h^m(x^-)-1 +
   \ind_f^\sigma(C^+)-\ind_f^\sigma(C^-) \cr
   &\ \ \ - \frac{\dim C^+}{2} + \frac{\dim C^-}{2} \cr
   &= \ind_{f,h}^\sigma(x^+)-\ind_{f,h}^\sigma(x^-)-1.
\end{align}
So we get the same formula for $\dim\MM(x^-,x^+;C_0,\dots,C_k)$ using
either Morse indices or signature indices. Since this dimension equals
zero for the moduli spaces contributing to the boundary operator in
Morse-Bott homology, this shows that the boundary operator has degree
$-1$ with respect to either grading. However, we
would like to point out that a mixture of Morse indices and
signature indices does not lead in general to a grading on Morse-Bott
homology unless all the connected components of the critical manifold
of $f$ have the same dimension. 
\\
Finally, consider the situation in Floer homology where the ambient
space is infinite-dimensional, but the components of
$\mathrm{crit}(f)$ and the moduli spaces $\MM(C_{i-1},C_i)$ are still
finite dimensional. Moreover, (under suitable hypotheses) the
dimension of these moduli spaces can be expressed by a formular
analogous to~\eqref{dim2} in terms of Conley-Zehnder indices:
$$
   \dim\MM(C_{i-1},C_i) = \mu_{CZ}(C_i)-\mu_{CZ}(C_{i-1}) +
   \frac{\dim C_i+\dim C_{i-1}}{2}.
$$
See e.g.~equation~\eqref{virdim2} for the Floer homology considered in this
paper. This suggest that the Conley-Zehnder index should be viewed as
a signature index, and the same computation as in the finited
dimensional case above yields the dimension formula
\begin{equation}\label{dim-casc3}
   \dim\MM(x^-,x^+;C_0,\dots,C_k) = \mu(x^+)-\mu(x^-)-1
\end{equation}
with respect to the signature index
$$
   \mu(x) := \mu_{CZ}(x)+\ind_h^\sigma(x). 
$$
Thus the boundary operator in Floer homology has degree $-1$ with
respect to $\mu$ and $\mu$ descends to an integer grading on Floer
homology. Actually, in the case considered in this paper this grading
takes values in $\frac{1}{2}+\Z$, where the shift by $\frac{1}{2}$
reflects the 1-dimensional constraint imposed on the free loop space.

\section[Spectral flow]{Spectral flow}\label{sec:spec}

We compare in this appendix two spectral flows which appear in Lagrange
multiplier type problems. To motivate this
we first consider the Lagrange multiplier functional in finite dimensions. 
Let $(M,g)$ be a Riemannian manifold and $(V,\langle \cdot,\cdot\rangle)$
be a Euclidean vector space. For functions $f \in C^\infty(M)$ and
$h \in C^\infty(M,V)$ the {\em Lagrange multiplier functional} 
$F \in C^\infty(M \times V)$ is given by
$$F(x,v)=f(x)+\langle v,h(x) \rangle.$$
For $v \in V$ we denote by $F_v \in C^\infty(M)$ the function given by
$$F_v=F(\cdot, v).$$
The Hessian of $F$ with respect to the metric 
$g \oplus g_V$ on $M \times V$, where $g_V=\langle\cdot,\cdot\rangle$
is the Euclidean scalar product on $V$, is given by
$$\Hess_F(x,v)=\left(\begin{array}{cc}
\Hess_{F_v}(x) & dh(x)^*\\
dh(x) & 0
\end{array}\right).$$
Here the adjoint of $dh(x)$ is taken with respect to the inner products
on $T_xM$ and $T_{h(x)}V \cong V$ given by the metric $g$ and by
$\langle \cdot, \cdot \rangle$.

We compare in this appendix the spectral flows of $\Hess_F$ and
$\Hess_{F_v}$ for Lagrange multiplier
functionals not necessarily defined on a finite dimensional
manifold. We will apply this in the following way.
For $F=\mathcal{A}^H$ the functional $F_v$ is
the action functional of classical mechanics whose spectral flow can be
computed via the Conley-Zehnder indices \cite{robbin-salamon}. 

To formulate our theorem we use the set-up of Robbin and Salamon
\cite{robbin-salamon}. 
Let $W$ and $H$ be separable real Hilbert spaces such that $W \subset
H$ is dense and the inclusion is compact. Let $A:W\to H$ be a bounded
linear operator (with respect to the norms on $W$ and $H$). Viewing
$A$ as an unbounded operator on $H$ with domain $\mathrm{dom}(A)=W$,
recall the following definitions (see e.g.~\cite{kato}). 
The adjoint operator 
$$
   A^*:\mathrm{dom}(A^*):=\{v\in H\mid |\langle v,Aw\rangle_H|\leq C|w|_H
   \text{ for all }w\in W\}\to H
$$
is defined by the equation 
$$
   \langle A^*v,w\rangle_H = \langle v,Aw\rangle_H.
$$
$A$ is called symmetric if $W\subset\mathrm{dom}(A^*)$ and
$A^*|_W=A$, i.e.~$\langle A^*v,w\rangle_H = \langle v,Aw\rangle_H$ for
all $v,w\in W$. $A$ is called self-adjoint if it is symmetric and
$\mathrm{dom}(A^*)=W$. The spectrum of $A$ is the set of all complex
numbers $\lambda$ such that the operator
$$A-\lambda \cdot\mathrm{id}\colon W \to H$$
is not bijective. Denote by $\mathrm{ker}(A)$ and $\mathcal{R}(A)$ the
kernel and range ($=$ image) of $A$, respectively. 
Denote by $\mathcal{L}(W,H)$ the space of bounded
linear operators and by 
$$
   \mathcal{S}(W,H)\subset\mathcal{S}(W,H)
$$
the subspace of self-adjoint operators. The following lemma clarifies
these concepts. 

\begin{lemma}\label{robsal}
Let $W \subset H$ be Hilbert spaces such that 
the inclusion is compact and let $A:W\to H$ be a symmetric bounded
linear operator. Then the following are equivalent:
\begin{description}
 \item[(i)] There exists $\lambda \in \mathbb{R}$ such that
  $A-\lambda \cdot \mathrm{id} \colon W \to H$
  is bijective.
 \item[(ii)] $A$ is self-adjoint considered as an unbounded operator
  on $H$ with domain $\mathrm{dom}(A)=W$.
 \item[(iii)] One of the defect indices 
  $d^\pm(A):=\mathrm{codim}(\mathcal{R}(A^\C\pm i\cdot 
  \mathrm{id}),H^\C)$ is zero. Here $A^\C:W^\C\to H^\C$ denotes the
  complex linear extension of $A$ to the complexified Hilbert spaces.   
 \item[(iv)] The spectrum of $A$ is discrete and consists of real
  eigenvalues of finite multiplicity.  
\end{description}
\end{lemma}

\textbf{Proof: }We first show that $(i) \Rightarrow (ii)$, i.e.~
$\mathrm{dom}(A^*)=\mathrm{dom}(A)=W$.
To see this let $v \in \mathrm{dom}(A^*)$. Since 
$A-\lambda \cdot \mathrm{id}$ is bijective and $A$ is symmetric, there
exists $w \in W$ such that
$$(A^*-\lambda \cdot \mathrm{id})v=(A-\lambda \cdot \mathrm{id})w
=(A^*-\lambda \cdot \mathrm{id})w,$$
which implies
$$(A-\lambda \cdot \mathrm{id})^*(v-w)=0.$$
Again using the fact that $A-\lambda \cdot \mathrm{id}$ is bijective,
we conclude that $(A-\lambda \cdot \mathrm{id})^*$ is bijective and
hence
$$v=w \in W.$$
It follows that $A$ is self-adjoint with $\mathrm{dom}(A)=W$. 
\\
If $A$ is self-adjoint, then both defect indices are zero
, see for example \cite[Theorem V.3.16]{kato}, so that
$(ii) \Rightarrow (iii)$.
\\
We show that $(iii) \Rightarrow (iv)$. Assume that $d^-(A)$ is zero, i.e.
$A^\C-i \cdot\mathrm{id}:W^\C\to H^\C$ is invertible. Since the inclusion
$\iota \colon W \to H$ is compact, the operator 
$$R:=\iota \circ(A-i \cdot \mathrm{id})^{-1} \colon H^\C \to H^\C$$
is compact. In particular, its spectrum $\sigma(R)$ 
consists of eigenvalues, the only accumulation point in $\sigma(R)$
is zero, and the eigenspace for each eigenvalue except zero is
finite dimensional. Let $\zeta \in \mathbb{C} \setminus \{i\}$. 
Then the following relations hold for the ranges
$$\mathcal{R}(A^\C-\zeta \cdot \mathrm{id})=R^{-1}
\mathcal{R}(R-\frac{1}{\zeta-i}\cdot\mathrm{id})$$
and the kernels
$$\mathrm{ker}(A^\C-\zeta \cdot \mathrm{id})=\mathrm{ker}
(R-\frac{1}{\zeta-i}\cdot \mathrm{id}).$$
In particular, we have a bijection
$$\sigma(R) \setminus \{0\} \to \sigma(A^\C), \quad
\mu \mapsto \frac{1}{\mu}+i$$
between the spectra under which the corresponding eigenspaces do not
change, i.e. for every $\mu \in \sigma(R)\setminus\{0\}$
the eigenspaces satisfy
$$E_\mu(R)=E_{\frac{1}{\mu}+i}(A^\C) \subset H.$$
We conclude that the spectrum of $A$ consists of discrete eigenvalues
of finite multiplicity, which are real because $A$ is symmetric. 
A similar argument holds for the case that $d^+(A)$ is zero.
This shows that $(iii)$ implies $(iv)$.
\\
That $(iv) \Rightarrow (i)$ is obvious. This finishes the proof
of the lemma. \hfill $\square$
\\ \\
Assume in addition that $V$ is a finite dimensional Hilbert space. Let
$A \in \mathcal{S}(W,H)$ be as before and 
$B \in \mathcal{L}(V,H)$ be a bounded linear operator. We denote
by $A_B \colon W \oplus V \to H \oplus V$ the bounded symmetric operator
defined by
$$A_B(w,v)=(Aw+Bv,B^*w).$$
As a consequence of Lemma~\ref{robsal} we get the following corollary.

\begin{cor}The operator $A_B$ is in $\mathcal{S}(W \oplus V,H\oplus V)$.
\end{cor}
\textbf{Proof: }By Lemma~\ref{robsal} we have to show that
$A_B$ is self-adjoint. This is true if $B=0$. For arbitrary $B$ 
this follows from Theorem V.4.3 in \cite{kato}. \hfill $\square$
\\ \\
In the following orthogonality is always understood with respect to the inner
product of $H$ and never with respect to the inner product of $W$.
\begin{fed}\label{regular}
Let $A \in \mathcal{S}(W,H)$, and $B \in \mathcal{L}(V,H)$.
We say that the tuple $(A,B)$ is {\em regular} if
\begin{description}
 \item[(i)] $B$ is injective.
 \item[(ii)]A maps $\mathcal{R}(B)$ to itself and
  the restriction $\hat{A}:=A|_{\mathcal{R}(B)}$ is bijective. 
\end{description}
\end{fed}
A regular pair $(A,B)$ gives rise to the symmetric form 
$$S_{A,B}:=B^*\hat{A}^{-1}B$$
on $V$ whose signature we denote by
\begin{equation}\label{sign}
\sigma(A,B)=\mathrm{sign}(S_{A,B}).
\end{equation}

We now consider maps $A \colon \mathbb{R} \to \mathcal{S}(W,H)$ and
$B \colon \mathbb{R} \to \mathcal{L}(V,H)$ which are continuous with
respect to the norm topology such that the limits
$$\lim_{s \to \pm \infty}A(s)=A^\pm,\quad
\lim_{s \to \pm \infty}B(s)=B^\pm$$
exist and $A^\pm \in \mathcal{S}(W,H)$. Accordingly we
define the map 
$A_B \colon \mathbb{R} \to \mathcal{S}(W \oplus V,H \oplus V)$ 
by
$$A_B(s):=A(s)_{B(s)}, \quad s \in \mathbb{R}.$$
Denote by 
$$\mathcal{A}=\mathcal{A}(\mathbb{R},W,H)$$
the space of maps $A \colon \mathbb{R} \to \mathcal{S}(W,H)$
as above, which in addition satisfy that $A^\pm$ is bijective.
Recall the following theorem of Robbin and Salamon about the existence
of the spectral flow \cite[Theorem 4.3]{robbin-salamon}.
\begin{thm}\label{specflow}
There exist unique maps $\mu \colon \mathcal{A}(\mathbb{R},W,H) \to 
\mathbb{Z}$, one for every compact dense injection of Hilbert spaces
$W \hookrightarrow H$, satisfying the following axioms.
\begin{description}
 \item[(homotopy)] $\mu$ is constant on connected components of
  $\mathcal{A}(\mathbb{R},W,H)$.
 \item[(constant)] If $A$ is constant, then $\mu(A)=0$.
 \item[(direct sum)] $\mu(A_1 \oplus A_2)=\mu(A_1)+\mu(A_2)$.
 \item[(normalization)] For $W=H=\mathbb{R}$ and $A(t)=\arctan(t)$,
  we have $\mu(A)=1$. 
\end{description}
The number $\mu(A)$ is called the \emph{\textbf{spectral flow}}
of $A$.
\end{thm}

These axioms easily imply the following generalization of the
(normalization) axiom:\\
{\bf (crossing) }{\em For $W=H$ finite dimensional,
$$
   \mu(A)=\frac{1}{2}\Bigl(\mathrm{sign}(A^+)-\mathrm{sign}(A^-)\Bigr).
$$}
To define the spectral flow also for
maps $A$ whose limits $A^\pm$ are not
necessarily bijective we choose a smooth cutoff function
$\beta \in C^\infty(\mathbb{R},[-1,1])$ such that $\beta(s)=1$ for
$s \geq 1$ and $\beta(s)=-1$ for $s \leq -1$ and define
\begin{equation}\label{extsf}
A_\delta:=A-\delta \beta \cdot \mathrm{id}, \quad 
\mu(A):=\lim_{\delta \searrow 0}\mu(A_\delta).
\end{equation}
Note that the limit in $\mu(A)$ stabilizes for sufficiently small
$\delta>0$. 

\begin{thm} \label{lagspec}
Let $A \colon \mathbb{R} \to \mathcal{S}(W,H)$ and
$B \colon \mathbb{R} \to \mathcal{L}(V,H)$ be continuous maps
whose limits $\lim_{s \to \pm \infty}A(s)=A^\pm$ and
$\lim_{s \to \pm \infty}B(s)=B^\pm$ exist. Assume moreover that 
$(A^\pm,B^\pm)$ are regular pairs. 
Then the spectral flows of $A$ and $A_B$ are related by
$$\mu(A_B)=\mu(A)+\frac{1}{2}\bigg(\sigma(A^-,B^-)-\sigma(A^+,B^+)\bigg).$$
\end{thm}
\textbf{Proof: }
Choose two cutoff functions 
$\beta^\pm \in C^\infty(\mathbb{R},[0,1])$ with the property that
$\beta^+(s)=1$ for $s \geq 1$, $\beta^+(s)=0$ for $s \leq 0$,
$\beta^-(s)=1$ for $s \leq -1$, and $\beta^-(s)=0$ for $s \geq 0$.
Define $S_{A,B} \in \mathcal{A}(\mathbb{R},V,V)$ by
$$S_{A,B}=\beta^+\cdot S_{A^+,B^+}+\beta^- \cdot S_{A^-,B^-}.$$
Abbreviate $P_V \colon H \oplus V \to V$ the canonical projection.
We prove the Theorem in three steps. 
\\ \\
\textbf{Step 1: }\emph{If  $\delta \neq 0$ is small
enough then $(A_B)_\delta$ is
homotopic to $(A_\delta)_B-\delta P_V^*S_{A,B}P_V$.}
\\ \\ 
To see this, note that $(A_B)_\delta=(A_\delta)_B-\delta P_V^*P_V$. So
it suffices to show that for $\delta>0$ sufficiently small and
symmetric linear maps $S^\pm \in \mathcal{L}(V)$ whose norm is small
enough the operators $(A^\pm_\delta)_{B^\pm}+P_V^*S^\pm P_V$ are
bijective.  
By Theorem V.4.3 in \cite{kato} the operators
$(A^\pm_\delta)_{B^\pm}+P_V^*S^\pm P_V$ are selfadjoint with dense
domain $W \oplus V$, hence by Lemma~\ref{robsal} their spectrum
consists of eigenvalues. Thus it suffices to show injectivity. 
Let $(w,v) \in W \times V$ be in the kernel of
$(A^\pm_\delta)_{B^\pm}+P_V^*S^\pm P_V$.
Then $(w,v)$ solves
\begin{equation}\label{as1}
\left.\begin{array}{c}
(A^\pm-\delta \cdot \mathrm{id})w+B^\pm v=0 \\
(B^\pm)^*w+S^\pm v=0
\end{array}
\right\}
\end{equation}
which is equivalent to
\begin{equation}\label{as2}
\left.\begin{array}{c}
w=-(A^\pm-\delta \cdot \mathrm{id})^{-1}B^\pm v\\
-(B^\pm)^*(A^\pm-\delta\cdot \mathrm{id})^{-1}B^\pm v
+S^\pm v=0.
\end{array}
\right\}
\end{equation}
But $(B^\pm)^*(A^\pm-\delta \cdot \mathrm{id})^{-1}B^\pm$ converges
to the nondegenerate linear map $S_{A^\pm,B^\pm}$ as $\delta$ goes to zero, 
and hence the second equation in (\ref{as2}) has only the trivial solution
$v=0$ if the norm of $S^\pm$ is small enough and hence $(w,v)=(0,0)$. 
This shows injectivity and hence the assertion of Step 1 follows. 
\\ \\
\textbf{Step 2:} \emph{For $\delta>0$ small
enough and $\epsilon \in [0,1]$ the maps 
$(A_\delta)_{\epsilon B}-\delta P_V^*S_{A,B}P_V$ are in
$\mathcal{A}(\mathbb{R},W \oplus V,H \oplus V)$, i.e.
their asymptotics are bijective.} 
\\ \\
Step 2 follows by a similar reasoning as Step 1. Assume that
$(w,v) \in W \oplus V$ lies in the kernel of one of the asymptotic operators. 
Then $(w,v)$ solves
$$\left.\begin{array}{c}
w=-\epsilon(A^\pm-\delta \cdot \mathrm{id})^{-1}B^\pm v\\
-\epsilon^2(B^\pm)^*(A^\pm-\delta\cdot \mathrm{id})^{-1}B^\pm v
-\delta S_{A^\pm,B^\pm} v=0
\end{array}
\right\}\;.
$$
Since both terms in the second equation have the same sign and the
first term converges to $-\epsilon^2S_{A^\pm,B^\pm}v$ as $\delta$ goes to
zero, these equations have only the trivial solution.
\\ \\
\textbf{Step 3: }\emph{We prove the theorem.}
\\ \\
Using the properties of the spectral flow
from Theorem~\ref{specflow} we are now in position to compute
\begin{eqnarray*}
\mu\big((A_B)_\delta\big)&=&\mu\big((A_\delta)_B-\delta P^*_VS_{A,B}P_V\big)\\
&=&\mu\big((A_\delta)_0-\delta P^*_VS_{A,B}P_V\big)\\
&=&\mu(A_\delta \oplus -\delta S_{A,B})\\
&=&\mu(A_\delta)+\mu(-\delta S_{A,B})\\
&=&\mu(A_\delta)+\frac{1}{2}\bigg(\sigma(A^-,B^-)-\sigma(A^+,B^+)\bigg).
\end{eqnarray*}
Here we have used Step 1 for the first equality, Step 2 for the second
one, and the (crossing) property of $\mu$ for the last one. 
Taking the limit $\delta \searrow 0$ the theorem follows. \hfill $\square$
\\ \\
There are scenarios where the signature $\sigma(A,B)$ can easily be computed.
We formulate such an example for a finite dimensional Lagrange multiplier
functional which can easily be generalized to infinite dimensional
examples.
\begin{lemma}\label{varlag}
Suppose that  $(M,g)$ is a Riemannian manifold, $f \in C^\infty(M)$, and
$h \in C^\infty(M)$ such that $0$ is a regular value of $h$. Let 
$(x_0,v_0)$ be a critical point of 
the Lagrange multiplier functional $F \in C^\infty(M \times \mathbb{R})$
given by $F(x,v)=f(x)+v \cdot h(x)$. Assume that
there exists $\epsilon>0$ and a smooth curve 
$(x,v) \in C^\infty\big((-\epsilon,\epsilon),M \times \mathbb{R}\big)$ 
satisfying 
$(x(0),v(0))=(x_0,v_0)$ such that the following holds
\begin{description}
 \item[(i)] $\partial_\rho x(0)=\nabla h(x_0)$,
 \item[(ii)] $\partial_\rho v(0) \neq 0$,
 \item[(iii)] $(x(\rho),v(\rho))$ for $\rho \in (-\epsilon, \epsilon)$
  is a critical point of
  the Lagrange multiplier functional $F^\rho \in C^\infty(M \times \mathbb{R})$
  given by $F^\rho(x,v):=f(x)+v\cdot(h(x)-\rho)$.
\end{description}
Then the pair
$(\Hess_{F_{v_0}}(x_0), \nabla h(x_0))$ is regular in the sense of
Definition~\ref{regular} and its signature is
$$\sigma\big(\Hess_{F_{v_0}}(x_0),\nabla h(x_0)\big)=
-\mathrm{sign}\big(\partial_\rho v(0)\big).$$
\end{lemma}
\textbf{Proof: }The identity
$$dF^\rho\big(x(\rho),v(\rho)\big)=0$$ 
for $\rho \in (-\epsilon, \epsilon)$ is equivalent to
\begin{equation}\label{rlag}
\left.\begin{array}{c}
df\big(x(\rho)\big)+v(\rho)\cdot dh\big(x(\rho)\big)=0\\
h\big(x(\rho)\big)=\rho.
\end{array}\right\}
\end{equation}
The first equation in (\ref{rlag}) can be written as
$$\nabla f\big(x(\rho)\big)+v(\rho)\cdot\nabla h\big(x(\rho)\big)=0.$$
Differentiating this identity with respect to $\rho$ and evaluating
at $\rho=0$ we compute
using assumption $(i)$
\begin{eqnarray*}
0&=&\Hess_f(x_0)\partial_\rho x(0)+v_0\Hess_h(x_0)\partial_\rho x(0)+
\partial_\rho v(0)\nabla h(x_0)\\
&=&\Hess_{F_{v_0}}(x_0)\nabla h(x_0)+\partial_\rho v(0) \nabla h(x_0).
\end{eqnarray*}
In particular, $\nabla h(x_0)$ is an eigenvector of 
$\Hess_{F_{v_0}}$ to the nonzero eigenvalue $-\partial_\rho v(0)$. 

It is now straightforward to check that the pair
$(\Hess_{F_{v_0}}(x_0), \nabla h(x_0))$ is regular. Condition $(i)$ in 
Definition~\ref{regular} follows from the assumption that
$0$ is a regular value of $h$ and thus $\nabla h(x_0)\neq 0$. Since
$\nabla h(x_0)$ is an eigenvector of the Hessian to a nonzero
eigenvalue, condition $(ii)$ is satisfied as well. 

To compute the signature we calculate
\begin{eqnarray*}
\sigma\big(\Hess_{F_{v_0}}(x_0), \nabla h(x_0)\big)&=&
\mathrm{sign} \big(dh(x_0) \hat{\Hess}_{F_{v_0}}(x_0)^{-1} \nabla h(x_0)\big)\\
&=&\mathrm{sign}\bigg(-\frac{||\nabla h(x_0)||^2}{\partial_\rho v(0)}\bigg)\\
&=&-\mathrm{sign}\big(\partial_\rho v(0)\big).
\end{eqnarray*}
This proves the lemma. \hfill $\square$

\section{Some topological obstructions}\label{app:top}

During the first author's talk at the Workshop on Symplectic Geometry,
Contact Geometry and Interactions in Lille 2007, E.~Giroux suggested
that in the case $n=2$, Corollary~\ref{cor:subcrit} results from
the following topological fact.

\begin{lemma}\label{n=2}
There exists no smooth embedding of
$S^*S^2\cong\R P^3$ into a subcritical Stein surface. 
\end{lemma}

Right after the talk, participants suggested the following three
proofs of this fact. Note that every subcritical Stein surface
is $\C^2$ or a boundary connected sum of copies of $S^1\times\R^3$, 
which embeds smoothly into $\R^4$; thus for the lemma it suffices to
prove that $\R P^3$ admits no embedding into $\R^4$.  

{\bf Proof 1 (V.~Kharlamov): }
This proof is based on the following theorem of Whitney (see
e.g.~\cite{massey}): The Euler number $e(\Sigma)\in\Z$ of the normal
bundle of a closed connected non-orientable surface $\Sigma$ embedded
in $\R^4$ satisfies $e(\Sigma)\equiv 2\chi(\Sigma)$ mod $4$. 
Now suppose we have an embedding $\R P^3\subset\R^4$. Then the normal
Euler number of the linear subspace $\R P^2\subset\R P^3\subset\R^4$
satisfies $e(\R P^2)\equiv 2$ mod $4$. But a nonvanishing normal
vector field to $\R P^3$ in $\R^4$ (which exists because $\R P^3$ is
orientable) provides a nonvanishing section of the normal bundle of
$\R P^2\subset\R^4$, contradicting nontriviality of $e(\R P^2)$.   
\hfill $\square$

{\bf Proof 2 (T.~Ekholm): }
This proof is based on the following theorem of Ekholm~\cite{ekholm}: 
The Euler characteristic of the (resolved) self-intersection surface
of a generic immersion of $S^3$ to $\R^4$ has the same parity as the
number of quadruple points. 
Now suppose we have an embedding $\R P^3\subset\R^4$. Composition with
the covering $S^3\to\R P^3$ yields an immersion of $S^3$ to $\R^4$
which can be perturbed (via a normal vector field to $\R P^3$
vanishing along $\R P^2$) to have self-intersection surface $\R P^2$
and no quadruple points, contradicting Ekholm's theorem. 
\hfill $\square$

The third proof, suggested by P.~Lisca, yields in fact the following
more general result: 

\begin{prop}\label{even}
For $n\geq 2$ even there exists no smooth embedding of 
$S^*S^n$ into a subcritical Stein $2n$-manifold. 
\end{prop}

{\bf Proof 3 (P.~Lisca): }
Suppose we have an embedding $S^*S^n\cong\Sigma\subset V$ into a
subcritical Stein $2n$-manifold $V$ for $n\geq 2$ even. From
$H_{2n-1}(V;\Z)=0$ it follows that $\Sigma$ bounds a compact subset 
$B\subset V$. Denote by $C$ the closure of $V\setminus B$,
so $V=B\cup_\Sigma C$. Since $n$ is even, the Gysin homology sequence
of the sphere bundle $S^*S^n\to S^n$ shows $H_n(S^*S^n;\Z)=0$ and
$H_{n-1}(S^*S^n;\Z)=\Z_2$. Using this and $H_n(V;\Z)=0$, the
Mayer-Vietoris sequence for $V=B\cup_\Sigma C$ implies
$H_n(B;\Z)=H_n(C;\Z)=0$. Let $X:=D^*S^n\cup_\Sigma B$ be the closed
oriented $2n$-manifold obtained by gluing the unit disk cotangent
bundle $D^*S^n$ (with orientation reversed) to $B$ along
$\Sigma$. Again by the Mayer-Vietoris sequence, we find that the free
part $H_n(X;\Z)/{\rm torsion}$ is isomorphic to $\Z$ and generated by
the zero section $S^n\subset D^*S^n$. But for $n$ even the zero
section $S^n$ has self-intersection number $2$ in $D^*S^n$, hence $-2$
in $X$, contradicting unimodularity of the intersection form (which is
an immediate consequence of Poincar\'e duality). 
\hfill $\square$

{\em Remark. }
The fact that $S^*S^2$ has no {\em exact contact} embedding into a
subcritical Stein surface $V$ can also be proved symplectically as
follows, using a deep result by Gromov about holomorphic
fillings. Since every subcritical Stein surface admits an {\em exact
symplectic} embedding into $\R^4$, it suffices again to consider the
case $V=\R^4$.  
Suppose there exists an exact contact embedding
$\iota:S^*S^2\hookrightarrow\R^4$. Removing the bounded component of
$\R^4\setminus\iota(S^*S^2)$ and gluing in the unit ball bundle $D^*S^2$
yields an exact convex symplectic manifold $W$ which contains an
embedded Lagrangian 2-sphere (the zero section in $D^*S^2$). On the
other hand, $W$ is symplectomorphic to $\R^4$ outside compact set. So
a result of Gromov~\cite{gromov} 
implies that $W$ is in fact symplectomorphic to $\R^4$. But this is a
contradiction because $\R^4$ does 
not admit any embedded Lagrangian 2-spheres.  

{\em Remark. }
We have not investigated obstructions to smooth embeddings of $S^*S^n$
into subcritical Stein manifolds for $n$ odd. As pointed out in the
introduction, at least for $n=3$ and $n=7$ there are no obstructions
and $S^*S^n$ embeds smoothly into $\C^n$.

\end{document}